\documentclass{article}
\usepackage{arxiv}

\usepackage{setspace}
\usepackage{lmodern}

\usepackage{textgreek}

\usepackage{fixltx2e}

\reversemarginpar

\usepackage{amsmath,amssymb,amsfonts,amsthm}
\usepackage{mathrsfs}

\usepackage[locale = FR]{siunitx} 
\DeclareSIUnit\vitesse{\meter\per\second}
\usepackage{eurosym}
\DeclareSIUnit{\octet}{o}

\usepackage{graphicx,array,tikz,multirow}
\usepackage{caption,subcaption} 
\usepackage{svg,float}
\usepackage{booktabs,paralist}
\newcolumntype{x}[1]{>{\centering\arraybackslash\hspace{0pt}}p{#1}}
\usepackage[section]{placeins}
\usepackage{hanging}

\usepackage{fancyhdr,emptypage} 

\fancypagestyle{plain}{ 
    \fancyhead{}\fancyfoot[C]{\thepage}

}

\usepackage[Conny]{fncychap}

\usepackage[francais,nohints,tight]{minitoc}		

\usepackage[nottoc]{tocbibind} 
\usepackage[titles]{tocloft}

\usepackage{textcomp} 

\usepackage{titling}
\usepackage{lipsum} 
\usepackage{csquotes} 

\usepackage{xspace}
\usepackage{afterpage}

\usepackage[textsize=footnotesize]{todonotes}

\usepackage{bookmark}
\usepackage{acronym}
\usepackage[nameinlink,french]{cleveref} 
\Crefname{figure}{Fig.}{Figs.} 
\crefname{figure}{fig.}{figs.}
\Crefname{equation}{Eq.}{Eqs.}
\crefname{equation}{eq.}{eqs.}
\Crefname{table}{Table.}{Tables.}
\crefname{table}{table.}{tables.}

\definecolor{color_ref}{rgb}{1.0, 0.13, 0.32} 
\definecolor{color_link}{rgb}{0.0, 0.0, 1.0}
\definecolor{curcolor}{rgb}{0.0, 0.0, 1.0} 
\definecolor{brightpink}{rgb}{1.0, 0.0, 0.5} 
\definecolor{navyblue}{rgb}{0.0, 0.0, 1.0}
\hypersetup{
	colorlinks=true, 
	pdfstartview=FitV, 
	urlcolor=brightpink, 
	linkcolor= navyblue, 
	citecolor=color_ref 
}

\usepackage{wrapfig} 
\usepackage[utf8]{inputenc} 
\usepackage[T1]{fontenc}    
\usepackage{hyperref}       
\usepackage{url}            
\usepackage{booktabs}       
\usepackage{amsfonts}       
\usepackage{nicefrac}       
\usepackage{microtype}      
\usepackage{lipsum}
\usepackage{graphicx}
\usepackage{xcolor}
\usepackage{booktabs}
\usepackage{siunitx}
 
\graphicspath{ {./images/} }
\usepackage{amsmath}
\usepackage{array}
\theoremstyle{definition}
\newtheorem{theorem}{Theorem}[section]

\newtheorem{algorithm}[theorem]{Algorithme}

\newtheorem{definition}[theorem]{Definition}

\newtheorem{notation}[theorem]{Notations}

\newtheorem{proposition}[theorem]{Proposition}
\newtheorem{remark}[theorem]{Remark}

\title{Optimizing Resource Allocation for COVID-19 Vaccination Planning within Long-Term Care Facilities: A Refined Multilevel Linear Programming Approach}
\author{
 Mustapha Kaci\\
 Department of Mathematics, Signal Image Parole (SIMPA) Laboratory\\ University of Oran Mohamed Boudiaf USTO-MB, Oran, Algeria.\\
  {\color{blue}\texttt{kaci.mustapha.95@gmail.com}}
   }

\usepackage{acronym}
\usepackage{tcolorbox}

\tcbuselibrary{theorems}
\newtcbtheorem{definitionbox}{Définition}{
  enhanced,
  colback=blue!10,
  colframe=blue!20,
  arc=5pt,
  boxrule=1pt,
  left=10pt,
  right=10pt,
  coltitle=black,
	number within=<section>,
}{def}
\newcommand{\Max}{\operatorname{Max}}
\newcommand{\sign}{\operatorname{sign}}

\usepackage{fancybox}
\usepackage{awesomebox}
\usepackage{mdframed}
\usepackage{algorithm}
\usepackage{comment}
\usepackage{algpseudocode}

\usepackage{multirow}
\usepackage{enumitem}

\usepackage{tikz}
\usetikzlibrary{positioning}
\usetikzlibrary{shapes, arrows}

\usepackage{tikz}
\usetikzlibrary{shapes,arrows,positioning}
\usepackage{sidecap, caption}
\usepackage{subcaption}

\usepackage{rotating}

\usepackage{booktabs}
\usepackage{diagbox}
\usepackage{multirow}
\begin{document}
\fontfamily{ptm}\selectfont
\maketitle
\begin{abstract}
In this paper, we introduce an algorithm designed to solve a Multilevel MOnoObjective Linear Programming Problem (ML(MO)OLPP). Our approach is a refined adaptation of Sinha and Sinha's linear programming method, incorporating the development of an "interval reduction map" that precisely refines decision variable intervals based on the influence of the preceding level's decision maker. Each construction stage is meticulously examined. The effectiveness of the algorithm is validated through a detailed numerical example, illustrating its practical applicability in resource management challenges.  With a specific focus on vaccination planning within long-term care facilities and its relevance to the COVID-19 pandemic, our study addresses the optimization of resource allocation, placing a strong emphasis on the equitable distribution of COVID-19 vaccines.
\end{abstract}


\tableofcontents

\section{Introduction}\label{section1}

Multilevel optimization, a robust method for modeling decision making scenarios with dynamic interactions among multiple decision makers within a hierarchical structure, has found widespread applications across various sectors \cite{kaci, kaci2023, abo}. In Multilevel mathematical programming problems, each level corresponds to a distinct decision maker, creating opportunities for cooperative optimization \cite{kaci, sinna2024, baky}. Emphasizing cooperative scenarios and ensuring a minimum level of satisfaction for lower-level decision makers are crucial for organizational benefit. Despite the prevailing emphasis on bilevel programming, the landscape of programming problems extends beyond this scope, prompting researchers to explore diverse approaches to address multilevel programming challenges \cite{l2, cite1, ref36, ref37}.

\subsection{Literature Review}

Over the past two decades, extensive exploration in the domain of multilevel programming/decision making has been documented through various survey papers \cite{BEN22, BEN39, BEN89, BEN126, BEN154}. These publications primarily focus on the early stages of research related to fundamental bilevel decision making, often within the confines of traditional solution approaches or specific application domains. The origins of multilevel decision making can be traced back to 1973 with the seminal paper by Bracken and McGill \cite{ref36, ref37}. Since the 1980s, a diverse body of related research has emerged, encompassing designations such as multilevel programming, multilevel optimization, and multilevel decision making.

Bilevel Programming (BLP) represents a significant category in mathematical programming, offering a structured approach to model and solve optimization problems involving two levels of decision making. Numerous research efforts have been dedicated to exploring the theoretical and applied aspects of BLP, revealing its features, diverse applications, and associated solution methods.

In various contexts, BLP has been studied for its capabilities in formulating piecewise linear functions and its connections to other optimization problems, as extensively detailed in \cite{BEN22}. This paper not only reviews prior findings but also introduces new results, addressing the NP-hardness of BLP and clarifying confusing representations in the literature. Another comprehensive introduction to BLP is provided by \cite{BEN39}, motivating this class through a simple application and detailing the general formulation of bilevel programs across linear, linear-quadratic, and nonlinear cases.

Special attention is given to applications in the energy networks domain, as outlined in \cite{BEN89}. New problems in this area, including natural gas cash-out, deregulated electricity market equilibrium, and biofuel design, are efficiently addressed as mixed-integer bilevel programs. 

The survey in \cite{BEN126} encompasses both theoretical and applied results, exploring emerging themes such as migration processes, allocation problems, information protection, and cybersecurity.
On the other hand, delves into decision making problems in decentralized organizations modeled as Stackelberg games, formulated as bilevel mathematical programming problems with two decision makers. The paper emphasizes the NP-hardness of solving such problems when decision makers lack motivation to cooperate. However, cooperative decision making is explored using interactive fuzzy programming, providing solutions aligned with decision makers' preferences. The paper extends these methods to bilevel linear programming problems in multi-objective environments and under uncertainty.

The bibliography compiled by Vicente and Calamai  \cite{BEN154} stands as a dynamic and permanent contribution to bilevel and multilevel programming references. The authors invite suggestions for additions, corrections, and modifications. The overview of past and current research in bilevel and multilevel programming provided in the summary facilitates researchers' understanding of the references in this area. As a valuable resource, this bibliography supports and encourages further research in multilevel programming and decision making.

\subsection{Recent Advances}

In recent years, the landscape of multilevel programming has witnessed the introduction of cutting-edge methodologies. Abo-Sinna's algorithm \cite{abo} targets the identification of $\alpha$-Pareto optimal solutions, incorporating tolerance membership functions and exploring cooperative scenarios within the BLP problems with fuzzy parameters. Simultaneously, the work in \cite{baky} contributes two innovative algorithms employing fuzzy goal programming for addressing MultiLevel MultiObjective Linear Programming Problem (ML-MOLPP), with a specific emphasis on optimizing fuzzy goals for decision makers. Additionally, the paper in \cite{baky2009} proposes a fuzzy goal programming algorithm to navigate the intricacies of solving decentralized bilevel multiobjective programming problems. Sinha and  Sinha \cite{sinha2002} utilize the fuzzy mathematical programming approach to minimize objectives using linear membership functions. Building on this foundation, their subsequent work in \cite{sinha} introduces a practical method for solving linear multilevel decentralized programming problems based on linear programming within a supervised search procedure.

Moreover, parallel efforts by Kaci and Radjef provides an intriguing perspective. In their paper \cite{kaci}, a new approach is developed to solve ML-MOLPPs. This approach, not considering the hierarchy of the problem, provides the set of all compromises. A simple criterion is established to test the feasibility of potential compromises, and a method based on Yu and Zeleny's approach is employed to generate the set of all compromises.  In a subsequent work \cite{kaci2023}, the authors extend their exploration by introducing a novel algorithm for ML-MOLPP. This algorithm, nested within the adaptive method of linear programming, begins by generating the set of all possible compromises (non-dominated solutions). Subsequently, it efficiently selects the best compromise among the potential settlements. Transforming the initial multilevel problem into an ML-MOLPP with bounded variables, the adaptive method is then applied, proving its efficiency compared to the simplex method. 

In addition to these advancements, two notable alternative techniques for solving ML-MOLPPs have been proposed. The first technique, presented in \cite{LM1}, introduces an approach based on Fuzzy Goal Programming (FGP), offering a simpler and computationally less intensive alternative to the algorithm proposed by Baky \cite{baky}. This technique transforms each objective function at each level into fuzzy goals, defining suitable membership functions and control vectors for each level's decision makers. The fussy goal programming approach is then employed to achieve the highest degree for each membership goal by minimizing the sum of negative deviational variables, ultimately providing a compromise optimal solution for ML-MOLPP problems. The second technique, detailed in \cite{LM2}, addresses fully neutrosophic multilevel multiobjective programming problems, employing a unique strategy to convert the problem into equivalent ML-MOLPPs with crisp values. The proposed approach utilizes ranking functions and fuzzy goal programming, demonstrating simplicity and uniqueness in providing compromise optimal solutions for decision makers. These alternative techniques contribute valuable options to the evolving landscape of methodologies for solving ML-MOLPPs.

Furthermore, the work in \cite{BAH} introduces a FGP algorithm for solving bilevel multiobjective programming problems with fuzzy demands. This algorithm addresses the experts' imprecise or fuzzy understandings of parameters in the problem formulation, assumed to be characterized as fuzzy numbers. The FGP algorithm, applied to achieve the highest degree for each membership goal, provides a satisfactory solution for all decision makers in the presence of fuzzy demands.

\subsection{Impact of COVID-19 and Associated Optimization Approaches}

The unparalleled challenges posed by the COVID-19 pandemic have reshaped decision making scenarios across various domains. The dynamic and often unpredictable nature of the crisis, from the surge in healthcare demands to disruptions in education and supply chains, necessitates a new level of adaptability. Traditional decision making frameworks are strained under the weight of uncertainty and urgency inherent in the pandemic's complexities. It is in this landscape that optimization emerges as a linchpin, offering a systematic and data-driven approach to navigate the intricacies of resource allocation, logistics, and strategic planning. The urgency of the pandemic underscores the critical role of optimization in ensuring efficient, equitable, and timely decisions, making it an indispensable tool for addressing the unique challenges precipitated by COVID-19.

A primary major concern is effectively managing hospital resources amid a spike in infections, while minimizing disruption to regular healthcare services. A multi-objective, multi-period linear programming model was proposed to optimize the distribution of infected patients, the evacuation rate of non-infected patients, and the creation of new COVID-19 intensive care units, simultaneously minimizing the total distance traveled by infected patients, the maximum evacuation rate of non-infected patients, and the infectious risk to healthcare professionals \cite{cov1}.

The education sector has also been profoundly impacted, with the need to adopt new learning approaches during school closures. A study was conducted to understand the learning habits of students in India (Maharashtra) during the suspension of classes due to the pandemic. This research used a linear programming model to analyze the importance of teachers during online learning \cite{cov2}.

In the context of the blood supply chain, a study sought to optimally configure a multi-echelon blood supply chain network considering uncertainties related to demand, capacity, and blood disposal rates. This approach used a bi-objective Mixed-Integer Linear Programming (MILP) model and interactive possibilistic programming to optimally address the problem considering the special conditions of the pandemic \cite{cov3}.

The pandemic has also had a significant impact on mental health, increasing the demand for psychological intervention services. A study combined statistical analyses and discrete optimization techniques to solve the problem of assigning patients to therapists for crisis intervention with a single tele-psychotherapy session. The integer programming model was validated with real-world data, and its results were applied in a volunteer program in Ecuador \cite{cov4}.

Finally, the equitable distribution of COVID-19 vaccines has become crucial for economic recovery, especially in developing countries. A mixed-integer linear programming model was proposed for equitable COVID-19 vaccine distribution, considering constraints such as specific refrigeration requirements for certain vaccines and the availability of storage capacity \cite{cov5}.

\subsection{Proposed Approach and Contributions}
This study introduces a novel algorithm tailored to solve ML(MO)OLPP. Drawing inspiration from the linear programming method developed by Sinha and Sinha, our approach offers a modified framework enriched by the adaptive method pioneered by Gabasov, Kirillova, and Kostyukova \cite{Gabasov,GABASOV1979, Gabasov2}. By incorporating the adaptive method's numerical constructive approach, our algorithm aims to enhance the efficiency and resolution speed compared to conventional methods such as the simplex method.

Building upon Sinha and Sinha's foundational work \cite{sinha}, our enhanced algorithm introduces a "interval reduction map" to refine decision variable intervals influenced by preceding decision makers in the multilevel structure. Notably, our methodology leverages the adaptive method to replace the simplex method, striving for accelerated problem-solving. The adaptive method, known for its effectiveness in solving monoobjective linear programming problems with bounded variables, introduces a numerical constructiveness that aligns seamlessly with the multilevel linear programming challenges presented in this work.

This paper presents a comprehensive algorithmic development, validated through a numerical demonstration showcasing its efficacy and practical applicability, notably in resource management contexts, with a specific focus on optimizing resource allocation and ensuring the equitable distribution of COVID-19 vaccines, particularly in the domain of vaccination planning within long-term care facilities. The adaptability of the algorithm to real-world challenges is underscored. Additionally, our study contributes to optimization methods, offering insights into the trade-offs between traditional linear programming and our adapted approach. The comparative analysis in Table~\ref{tab6} sheds light on computational performance variations under different parameter configurations ($\alpha_{1}$ and $\alpha_{2}$), aiding in the selection of optimization approaches based on specific requirements and resource constraints. With a focus on addressing resource management challenges, especially amid the ongoing COVID-19 pandemic.

\subsection{Paper Organization}

The organization of this paper is structured as follows: Section \ref{section2}, titled "Preliminaries," lays the foundational groundwork by presenting essential background information, setting the stage for a comprehensive understanding of subsequent developments. In Section \ref{section3}, we meticulously formulate the ML(MO)OLPP mathematically, offering a detailed exploration of its structural components. Progressing further, Section \ref{section4} narrows its focus to the linear formulation specific to the $p$th level, unraveling the intricacies of our proposed approach. Section \ref{section5} delves into the details of an algorithmic solution for solving the ML(MO)OLPP, utilizing Gabasov's Adaptive Method. Shifting the spotlight to practical application, Section \ref{section6} showcases our approach in the context of vaccination planning within Long-Term Care Facilities, illustrated through a carefully chosen numerical example. Finally, Section \ref{section7}, our "Conclusion," brings together our contributions, summarizes key findings, and explores the implications of our work, providing a comprehensive conclusion to our exploration of the ML(MO)OLPP and its real-world applications.

\section{Preliminaries}\label{section2}
\subsection{The Gabasov's adaptive method}

Gabasov's adaptive method relies on the linear search technique to find an optimal distance within the feasible region $S$ along a given direction. This approach, commonly employed in constrained optimization algorithms, involves several key steps. First, a direction $d \in \mathbb{R}^{n}$ is chosen to improve the objective function's value. Subsequently, the maximum distance that can be traveled along this direction without exiting the feasible region $S$ is determined. A linear search is then conducted in the chosen direction $d$ using the calculated maximum distance. The feasible solution $x$ is updated by moving optimally along the direction $d$. These steps are repeated until a specified termination criterion is met.

Consider a linear programming problem with bounded variables, formulated in the canonical form as follows:
\begin{equation}\label{Pbm01}
\left\{
\begin{array}{l}
      \underset{x}{\max} \hspace{0.12cm} f(x)=c^{T}x           \\
      Ax= b\\
			l\leq x\leq u,
\end{array}
\right.
\end{equation}
where $c$, $l$, $u$ $\in\mathbb{R}^{n}$, such that $\left\|l\right\| < \infty$ and $\left\|u\right\| < \infty$, and $b\in\mathbb{R}^{m}$. The matrix $A$ has dimensions $m \times n$ with $m < n$. 
\begin{definition}$\;$\\
\begin{enumerate}
\item A vector $x$ satisfying the constraints of problem (\ref{Pbm01}) is termed a {\it feasible solution} or simply a {\it plan}. The set of all feasible solutions is the {\it feasible region}, defined as:
$$
S:=\left\{x\in\mathbb{R}^{n}\hspace{0.2cm}:\hspace{0.2cm}Ax=b,\hspace{0.2cm}l\leq x\leq u\right\}.
$$

\item Let $J_{B}$ be a set of indices such that the cardinality of $J_{B}$ is equal to $m$, then $J_{B}$ is termed {\it support} if the matrix $A_{B}$ is invertible.

\item The pair $\left\{x, J_{B}\right\}$, consisting of a plan $x$ and a support $J_{B}$, is a {\it supporting plan}. A supporting plan is {\it non-degenerate} if:
$
l_{j}<x_{j}<u_{j},\hspace{0.25cm}j\in J_{B}.
$
\item A feasible solution \(x^{0}\) is {\it optimal} if \(f(x^{0}) = \underset{x\in S}{\max} \hspace{0.08cm} f(x)\), while \(x_{\epsilon}\) is {\it $\epsilon$-optimal} or {\it sub-optimal} if \(f(x^{0}) - f(x_{\epsilon}) \leq \epsilon\), where \(x^{0}\) is an optimal solution for problem (\ref{Pbm01}), and \(\epsilon\geq0\). 

\end{enumerate}
\end{definition}
The principle of Gabasov's adaptive method is as follows: start with an initial supporting plan $\left\{x, J_{B}\right\}$, which is a candidate solution satisfying the problem's constraints. Then, iterate to reach a new supporting plan $\left\{\bar{x}, \bar{J}_{B}\right\}$ while improving the objective function value, i.e., by seeking to satisfy the inequality $f(\bar{x}) \geq f(x)$.

\subsubsection{Suboptimality Estimate}
Consider $\left\{x, J_{B}\right\}$ as a supporting plan for problem (\ref{Pbm01}), and let $x^{0}$ be any plan such that $x^{0}=x+\Delta x$. The increase in the objective function is expressed as:
\begin{equation}\label{accf}
\Delta f(x)=f(x^{0})-f\left(x\right)=c^{T}x^{0}-c^{T}x=c^{T}\left(x^{0}-x\right)=c^{T}\Delta x.
\end{equation}
Furthermore, observe that:
$
A\Delta x=A\left(x^{0}-x\right)=Ax^{0}-Ax=b-b=0.
$
This implies:
$
A_{B}\Delta x_{B}+A_{N}\Delta x_{N}=0.
$
Consequently, $\Delta x_{B}$ can be expressed in terms of $\Delta x_{N}$:
$
\Delta x_{B}=-A_{B}^{-1}A_{N}\Delta x_{N}.
$
By substituting $\Delta x_{B}$ into equation (\ref{accf}), we obtain:
$
\Delta f(x)=c^{T}\Delta x=c^{T}_{B}\Delta x_{B}+c^{T}_{N}\Delta x_{N}
					                          =-c^{T}_{B}A_{B}^{-1}A_{N}\Delta x_{N}+c^{T}_{N}\Delta x_{N}.
$
This allows us to rewrite the increase in the objective function as:
\begin{equation}\label{accf2}
\Delta f(x)=-\left(c^{T}_{B}A_{B}^{-1}A_{N}-c^{T}_{N}\right)\Delta x_{N}.
\end{equation}

To simplify our analysis further, we denote the vectors $K$ and $E$ as follows:
$K^{T}_{B}=c_{B}^{T}A_{B}^{-1}$, $E^{T}=K^{T}A-c^{T}=(E^{T}_{N},E^{T}_{B})$,
where 
 $E_{N}^{T}=c_{B}^{T}A_{B}^{-1}A_{N}-c^{T}_{N}$, 
 $E_{B}^{T}=c_{B}^{T}A_{B}^{-1}A_{B}-c^{T}_{B}=0$.

Finally, using these notations, formula (\ref{accf2}) can be rewritten as follows:
\begin{equation}\label{accf3}
\Delta f(x)=-E^{T}\Delta x_{N}=\sum_{j\in J_{N}}E_{j}\left(x_{j}-x_{j}^{0}\right).
\end{equation}

\begin{theorem}[Optimality Criterion \cite{Gabasov,GABASOV1979}] \label{opcri} Let $\left\{x, J_{B}\right\}$ be a supporting plan of the problem $(\ref{Pbm01})$. Then, the following relations are sufficient to ensure the optimality of the supporting plan $\left\{x, J_{B}\right\}$:
\begin{equation}\label{CD2023}
\left\{\begin{array}{cc}
\begin{array}{lll}
E_{j}\geq0& \text{if} &x_{j}=l_{j} \\
E_{j}\leq0& \text{if} & x_{j}=u_{j} \\
E_{j}=0   & \text{if} & l_{j}<x_{j}<u_{j} 
\end{array},& j\in J_{N}.
\end{array}\right.
\end{equation}
Moreover, these relations are also necessary when the supporting plan $\left\{x, J_{B}\right\}$ is non-degenerate.
\end{theorem}

Using Theorem \ref{opcri}, we have the following inequality:
$x_{j}-u_{j}\leq x_{j}-x^{0}_{j}\leq x_{j}-l_{j}$, $j\in J_{N}$.
This allows us to establish the following upper bounds:
$$
\left\{\begin{array}{ccc}
E_{j}\left(x_{j}-x^{0}_{j}\right)\leq E_{j}\left(x_{j}-l_{j}\right),&\text{if}&E_{j}>0\vspace{0.15cm}\\
E_{j}\left(x_{j}-x^{0}_{j}\right)\leq E_{j}\left(x_{j}-u_{j}\right),&\text{if}&E_{j}<0.
\end{array}\right.
$$
Consequently, we derive the following upper bound, denoted as  $\beta(x,J_{B})$, termed as the suboptimality estimate:
$$
\begin{array}{cclcl}
f(x^{0})-f(x)=\sum_{j\in J_{N}}E_{j}\left(x_{j}-x_{j}^{0}\right)             &=&\sum_{E_{j}>0, j\in J_{N}}{E_{j}\left(x_{j}-x_{j}^{0}\right)}&+&\sum_{E_{j}<0, j\in J_{N}}{E_{j}\left(x_{j}-x_{j}^{0}\right)} \vspace{0.25cm}\\
             &\leq&\sum_{E_{j}>0, j\in J_{N}}{E_{j}\left(x_{j}-l_{j}\right)}&+&\sum_{E_{j}<0, j\in J_{N}}{E_{j}\left(x_{j}-u_{j}\right)}=\beta(x,J_{B}).
\end{array}
$$


\begin{theorem}[Sufficient Optimality Condition \cite{Gabasov,GABASOV1979}]
Let $\left\{x, J_{B}\right\}$ be a supporting plan of the problem (\ref{Pbm01}), and $\epsilon$ be a strictly positive number. Then, the plan $x$ is an $\epsilon$-optimal solution if and only if $\beta(x, J_{B}) < \epsilon$.
\end{theorem}

\begin{algorithm}[H]
\caption{Gabasov's monocriteria adaptive method}
\label{ama}

\textbf{Input:} Initial supporting plan $\left\{x, J_{B}\right\}$, parameter $\epsilon$.\\
\textbf{Output:} Optimal supporting plan $\left\{x, J_{B}\right\}$ or $\epsilon$-optimal plan.

\bigskip

{\textbf{Step 1:}} Compute $K^{T}_{B} = c_{B}^{T}A_{B}^{-1}$, $E^{T}_{j} = K^{T}a_{j} - c_{j}^{T}$, $j \in J_{N}$.\\
{\textbf{Step 2:}} Compute $\beta\left(x,J_{B}\right)$ using formula (\ref{formule0}).\\
{\textbf{Step 3:}} If $\beta\left(x,J_{B}\right) = 0$, stop with the optimal supporting plan {$\left\{x, J_{B}\right\}$}.\\
{\textbf{Step 4:}} If $\beta\left(x,J_{B}\right) \leq \epsilon$, stop with the $\epsilon$-optimal supporting plan {$\left\{x, J_{B}\right\}$}.\\
{\textbf{Step 5:}} Determine the set of non-optimal indices: 
\begin{center}
$J_{NNO} = \left\{j \in J_{N}: \left[E_{j} < 0, x_{j} < u_{j}\right] \text{ or } \left[E_{j} > 0, x_{j} > l_{j}\right]\right\}$.
\end{center}
{\textbf{Step 6:}} Choose an index $j_{0}$ from $J_{NNO}$ such that $\left|E_{j_{0}}\right| = \underset{j \in J_{NNO}}{\max}\left|E_{j_{0}}\right|$.\\
{\textbf{Step 7:}} Compute the admissible improvement direction $d$ using the relations:
\begin{equation}\label{formule2}
\left\{\begin{array}{lcc}
d_{j_{0}} = -\text{sign}\left(E_{j_{0}}\right) & &\vspace{0.15cm}\\
d_{j} = 0 & \text{if} & j \neq j_{0},\hspace{0.1cm}j \in J_{N}.\vspace{0.15cm}\\
d_{B} = -A_{B}^{-1}A_{j_{0}}d_{j_{0}} &&
\end{array}\right.
\end{equation}
{\textbf{Step 8:}} Compute  $\theta_{j_{1}} = \underset{j \in J_{B}}{\min}\hspace{0.1cm}\theta_{j}$, where $\theta_{j}$ is determined by the formula:
\begin{equation}\label{formule3}
\theta_{j} = \left\{\begin{array}{lll}
                                 \frac{u_{j} - x_{j}}{d_{j}},   &\text{if} & d_{j} > 0\vspace{0.15cm}\\
																 \frac{l_{j} - x_{j}}{d_{j}},   &\text{if} & d_{j} < 0\vspace{0.15cm}\\
																       \infty,                    &\text{if} & d_{j} = 0.
                       \end{array}   \right.
         \end{equation}.
{\textbf{Step 9:}} Compute $\theta_{j_{0}}$ using the formula:	
\begin{equation}\label{formule4}
\theta_{j_{0}} = \left\{\begin{array}{lll}
                                    x_{j_{0}} - l_{j_{0}} & \text{if} & E_{j_{0}} > 0\vspace{0.15cm}\\
																   u_{j_{0}} - x_{j_{0}} & \text{if} & E_{j_{0}} < 0.
																
                       \end{array}\right.                       \end{equation}	
{\textbf{Step 10:}} Compute $\theta^{0} = \min\left\{\theta_{j_{1}},\theta_{j_{0}}\right\}$.\\
{\textbf{Step 11:}} Compute $\overline{x} = x + \theta^{0}d$.\\
{\textbf{Step 12:}} Compute $\beta\left(\overline{x},J_{B}\right) = \beta\left(x,J_{B}\right) - \theta^{0}\left|E_{j_{0}}\right|$.\\
{\textbf{Step 13:}} If $\beta\left(\overline{x},J_{B}\right) = 0$, stop with the optimal supporting plan {$\left\{x, J_{B}\right\}$}.\\
{\textbf{Step 14:}} If $\beta\left(\overline{x},J_{B}\right) \leq \epsilon$,  stop with the supporting plan {$\left\{x, J_{B}\right\}$} $\epsilon-$optimal.\\
{\textbf{Step 15:}} If $\theta^{0} = \theta_{j_{0}}$, then set $\overline{J}_{B} = J_{B}$.\\
{\textbf{Step 16:}} If $\theta^{0} = \theta_{j_{1}}$, then set $\overline{J}_{B} = \left(J_{B} \backslash \left\{j_{1}\right\}\right)\cup\left\{j_{0}\right\}$.\\
{\textbf{Step 17:}} Set $x = \overline{x}$ and $J_{B} = \overline{J}_{B}$, and go to \textbf{Step 1}.
\end{algorithm}

\subsection{Multilevel Model for Efficient COVID-19 Vaccine Distribution}
We address a multilevel  allocation problem for the efficient management of resources during the COVID-19 pandemic. The problem is formulated using a three-level linear programming model, taking into account central, regional, and local management.

At the first level, central management, represented by the government of the United Kingdom, aims to maximize regional coordination and vaccine supply while minimizing overall costs. Possible decisions at this level include the quantity of vaccine doses allocated to the United Kingdom. Associated constraints relate to the national vaccine production capacity.

At the second level, regional management, each region of the United Kingdom manages resource allocation based on specific needs and capacities. Possible decisions at this level include the number of hospital beds allocated to each region and the number of resources allocated to each hospital in the region. The objective is to maximize vaccine utilization by region. Constraints include regional hospital capacity, regional resource stock, bed needs by region, and resource needs by hospital in each region.

At the third level, local management, each hospital within each region manages resource allocation to ensure adequate medical care for patients. Possible decisions at this level include the number of hospital beds allocated to each care unit in each hospital in the region. The objective is to maximize vaccine utilization by hospitals in the region. Constraints include hospital bed capacity for each care unit in each hospital, resource stock for each hospital in each region, and bed needs for each care unit.

This multilevel modeling allows for hierarchical resource management, with specific objectives at each management level. The associated figure illustrates the hierarchy and relationships between management levels, showing how decisions and objectives at the central level impact the regional and local levels. This approach provides a structured framework for the optimal allocation of resources during the COVID-19 pandemic.

The application of the developed method on a specific case featuring linear constraints and objective functions will be meticulously detailed in Section \ref{MDLF}. Subsequently, a practical numeric demonstration will be provided in Section \ref{MDLFA}, utilizing real-world data.

\begin{figure}
    \centering

    \tikzstyle{level} = [rectangle, rounded corners, line width=2pt,text centered,text width=3cm, draw=black,, fill=white]
    \tikzstyle{decision} = [rectangle, rounded corners,line width=2pt, text centered,text width=4cm, draw=black, fill=white]
    \tikzstyle{block} = [rectangle, draw, text width=3cm,line width=2pt, text centered, fill=white]
    \tikzstyle{line} = [draw, -latex]

    \begin{tikzpicture}[centered, node distance=1.65cm]
        \node [level] (level1) {{\bf Level 1:} Central Management (Government of the United Kingdom)};

        \node [level, below=of level1] (level2) {{\bf Level 2:} Regional Management (Regions of the United Kingdom)};
				
        \node [level, below=of level2] (level3) {{\bf Level 3:} Local Management (Hospitals within Each Region)};
				
        \node [decision, left=of level2] (decision1) {{\bf Objective:} Maximize coordination and vaccine supply while minimizing costs.\\Constraints: Budget constraints};
				
        \node [decision, right=of level2] (decision2) {{\bf Objective:} Maximize vaccine utilization by region.\\Constraints: Regional hospital capacity, resource stock, bed needs, and resource needs by hospital.};

        \node [decision, below=of decision1] (decision3) {{\bf Objective:} Maximize vaccine utilization by hospitals in the region.\\Constraints: Hospital bed capacity, resource stock, and bed needs for each care unit.};

        \path [line] (level1) -| (decision1);
        \path [line] (level2) -- (decision2);
        \path [line] (level3) -- (decision3);
				\path [line] (decision2) |- (level3);
        \path [line] (decision1) -- (level3);
				\draw [line, bend right] (decision3) to (level3);
    \end{tikzpicture}
    
    \caption{Detailed Multilevel Model for Equitable COVID-19 Vaccine Distribution in the United Kingdom\label{dim}}
\end{figure}
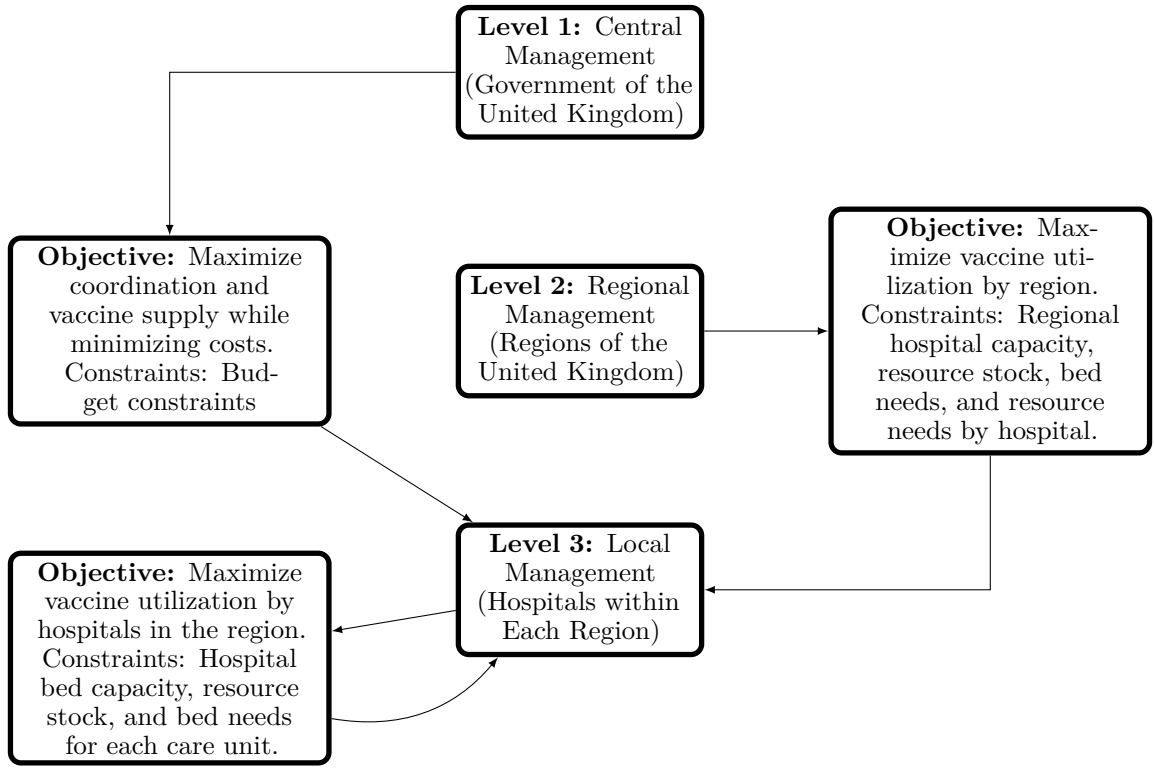

\begin{remark}
     In the model depicted in figure \ref{dim},  the arrows indicate the relationships between the levels of management. Decisions and objectives at the central level (level 1) impact the regional (level 2) and local (level 3) levels. Resources and constraints propagate from the central level to the regional and local levels, where more detailed decisions are made to optimize vaccine distribution.
\end{remark}

\section{Formulation of a ML(MO)OLPP}\label{section3}

Let's  denote the $p$th decision maker by $DM_{p}$, where $p=1\ldots P$ and $P\geq2$. These decision makers ($DM$s) are tasked with controlling the decision variables $\overline{x}^{p} = x_{p,1}, \ldots, x_{p,n_{p}} \in \mathbb{R}^{n_{p}}$, for $j=1\ldots n_{p}$, and $n = n_{1} + \ldots+ n_{P}$.
The decision variable column vector $x$ is defined as $\overline{x}^{1}, \ldots, \overline{x}^{P}$.\\
Consider $c_{p}^{i,j}$ as constants for $i, p = 1, \ldots, P$ and $j = 1, \ldots, n_{i}$. The objective function $f_{p} : \mathbb{R}^{n_{1}} \times \ldots \times \mathbb{R}^{n_{P}} \longrightarrow \mathbb{R}$, corresponding to the $p$-th level of our ML(MO)OLPP, is expressed as:
$$
\begin{array}{lclclcl}
f_{p}(x) &=& c_{p}^{T}x &=& \left(c_{p}^{1,j}\, c_{p}^{2,j}\, \, \, \ldots \, \, \, c_{p}^{P,j}\right) x&=&\left(c_{p}^{1,1}\, \, \, \ldots \, \, \, c_{p}^{1,n_{1}}\, c_{p}^{2,1}\, \, \, \ldots \, \, \,  c_{p}^{2,n_{2}}\, \, \, \ldots \, \, \, c_{p}^{P,n_{P}}\right) x.
\end{array}
$$
 This can be further expanded as:
$$
\begin{array}{lcl}
f_{p}(x) &=&c_{p}^{1,1}x_{1,1}+\ldots+c_{p}^{1,n_{1}}x_{1,n_{1}}+c_{p}^{2,1}x_{2,1}+\ldots+c_{p}^{2,n_{2}}x_{2,n_{2}}+\ldots+\ldots+c_{p}^{P,1}x_{P,1}+\ldots+c_{p}^{P,n_{P}}x_{P,n_{P}}.
\end{array}
$$
A $P$-level linear programming problem with a single objective function at each level can be formulated as follows:

\begin{equation}\label{mlpp1js}
\begin{array}{cl}
\textbf{Level 1:} 
& {\color{blue}\underset{\color{orange}\color{orange}\overline{x}^{1}}{\Max}}\hspace{0.1cm} f_{1}(x)={\color{blue}\underset{\color{orange}\overline{x}^{1}}{\Max}}\hspace{0.1cm}f_{1}({\color{orange}\overline{x}^{1}},\overline{x}^{2},\ldots,\overline{x}^{P})={\color{blue}\underset{\color{orange}\overline{x}^{1}}{\Max}}\hspace{0.1cm}c^{T}_{1}x, \\
& \text{such that } \overline{x}^{2},\ldots,\overline{x}^{P} \text{solves: }\vspace{0.2cm}\\
\textbf{Level 2:} 
& {\color{blue}\underset{\color{orange}\overline{x}^{2}}{\Max}}\hspace{0.1cm} f_{2}(x)={\color{blue}\underset{\color{orange}\overline{x}^{2}}{\Max}}\hspace{0.1cm}f_{2}(\overline{x}^{1},{\color{orange}\overline{x}^{2}},\ldots,\overline{x}^{P})={\color{blue}\underset{\color{orange}\overline{x}^{2}}{\Max}}\hspace{0.1cm}c^{T}_{2}x,\\
\vdots& \hspace{0.85cm}\vdots\hspace{2.9cm}\vdots\hspace{3.15cm}\vdots\\
\vdots& \hspace{0.85cm}\vdots\hspace{2.9cm}\vdots\hspace{3.15cm}\vdots\\
\vdots& \text{such that } \overline{x}^{P} \text{solves:}\hspace{0.87cm}\vdots\hspace{3.15cm}\vdots\vspace{0.2cm}\\
\textbf{Level P:} 
& {\color{blue}\underset{\color{orange}\overline{x}^{P}}{\Max}}\hspace{0.1cm} f_{P}(x)={\color{blue}\underset{\color{orange}\overline{x}^{P}}{\Max}}\hspace{0.1cm}f_{P}(\overline{x}^{1},\overline{x}^{2},\ldots,{\color{orange}\overline{x}^{P}})={\color{blue}\underset{\color{orange}\overline{x}^{P}}{\Max}}\hspace{0.1cm}c^{T}_{P}x,
\end{array}
\end{equation}

subject to:
$$
x\in S:=\left\{x\in\mathbb{R}^{n}\hspace{0.15cm}:\hspace{0.15cm}  Ax\leq b, x\geq 0, b\in\mathbb{R}^{m}\right\},
$$
where $S \neq \emptyset$ is the convex feasible region, $m$ is the number of constraints, $A$ is an $m \times n$ matrix, and $b$ is an $m$-dimensional vector.

\subsection{Resolution approach for multilevel problem (\ref{mlpp1js})}
Multilevel programming problems are marked by a hierarchical structure of $DM$s, where each decision level considers the decisions of the previous level to solve its own multiobjective problem. The following outlines a practical approach to solving such problems:
\begin{enumerate}
\item 
\textbf{Definition of decision levels:}
The first level (upper level) is the upper $DM$ ($DM_{1}$), and subsequent levels are subordinate $DM$s ($DM_{2}$, $\ldots$ , $DM_{P}$). For this type of problem, the variables $\overline{x}^{1}$, $\overline{x}^{2}$$\ldots$$\overline{x}^{P}$ must adhere to global constraints that define the feasible region $S$ of the multilevel problem. 
\item 
\textbf{Resolution of the upper $DM$'s problem (upper level):}
The $DM_{1}$ solves its monoobjective problem considering only the initial constraints. Denoting $f_1$ as the objective function of the 1$st$ level, the problem can be formulated as follows:
\begin{equation}
\begin{aligned}
{\color{blue}\underset{\color{orange}\color{orange}\overline{x}^{1}}{\Max}}\hspace{0.1cm} f_{1}(x)={\color{blue}\underset{\color{orange}\overline{x}^{1}}{\Max}}\hspace{0.1cm}f_{1}({\color{orange}\overline{x}^{1}},\overline{x}^{2},\ldots,\overline{x}^{P})={\color{blue}\underset{\color{orange}\overline{x}^{1}}{\Max}}\hspace{0.1cm}c^{T}_{1}x,\hspace{0.25cm} x\in S,
\end{aligned}
\end{equation}
where $\overline{x}^{1}$ are the $DM_{1}$'s decision variables, and $\overline{x}^{2}, \ldots, \overline{x}^{P}$ represents the decision variables of subsequent levels.
\item 
\textbf{Implementation of additional constraints in subordinate levels:}
After attaining its optimal solution $\overset{c}{\overline{x}^{1}}$, the upper level utilizes this solution to introduce additional constraints on the subordinate levels. These constraints narrow down the search space to a subset $S_{\bf 1} \subset S$ of feasible solutions.
\item 
\textbf{Resolution of subordinate $DM$s' problems}
Subordinate $DM$s then solve their own monoobjective problems while considering the new constraints imposed by the upper levels. The $DM_{p}$ solves their problem with $S_{\bf 1}$ as the search space for their decision variables. The problem for $DM_{p}$ can be formulated as follows:
\begin{equation}
\begin{aligned}
{\color{blue}\underset{\color{orange}\color{orange}\overline{x}^{p}}{\Max}}\hspace{0.1cm} f_{p}(x)={\color{blue}\underset{\color{orange}\overline{x}^{p}}{\Max}}\hspace{0.1cm}f_{p}(\overline{x}^{1},\ldots,{\color{orange}\overline{x}^{p}},\ldots,\overline{x}^{P})={\color{blue}\underset{\color{orange}\overline{x}^{p}}{\Max}}\hspace{0.1cm}c^{T}_{p}x,\hspace{0.25cm} x\in S_{\bf 1}.
\end{aligned}
\end{equation}

\item 
\textbf{Iteration of the Process:}
The process continues by iterating steps 3 and 4 for each subordinate level until the last level is reached.
\item 
\textbf{Attainment of Compromise:}
Once all levels have solved their problems, progressively improving the solution from the upper level. The final compromise is reached at the last level, where all decisions are considered in a balanced manner, taking into account the constraints imposed by higher levels. Thus, the final compromise is the solution at the last level, reflecting the preferences of each decision level and resulting in an effective and balanced solution for ML(MO)OLPP.
\end{enumerate}

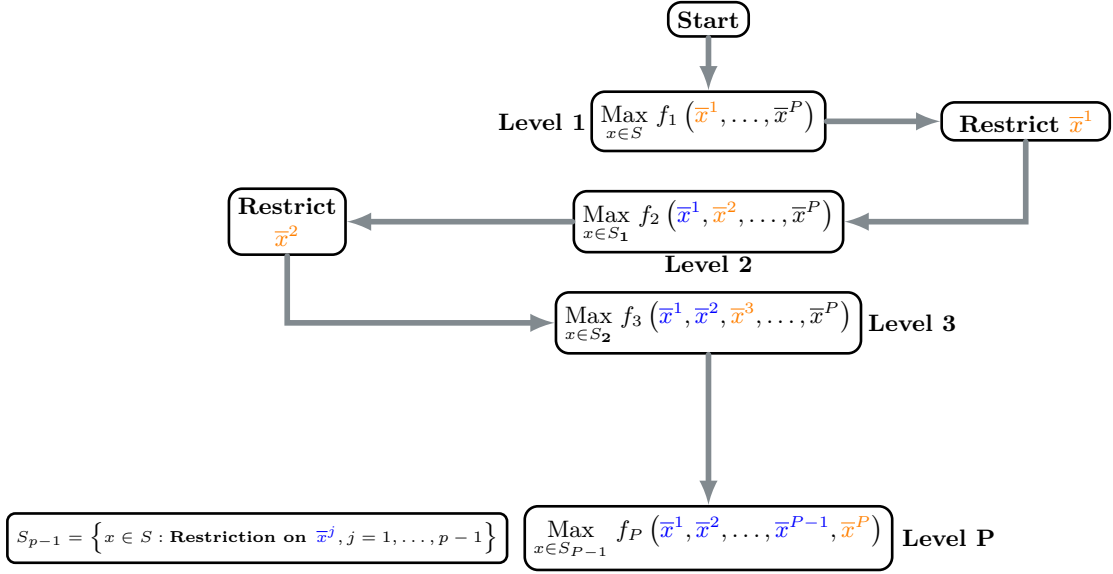
\begin{figure}[H]

\begin{tikzpicture}
    \tikzset{
      node distance=5mm and 0mm,
      boite/.style={
        fill=blue!70!cyan!60,draw,align=center,
        font=\footnotesize\bfseries,text=white, 
      },
			 boite4/.style={
        fill=white,draw,align=center,
        font=\footnotesize\bfseries,line width=1pt,text=black, 
      },
			boite2/.style={
        fill=white,draw,align=center,
        font=\footnotesize\bfseries,line width=1pt,text=black
      },
			boite3/.style={
        fill=green!70!cyan!60,draw,align=left,
        font=\footnotesize\bfseries,text=black, minimum width=2cm, minimum height=1cm
      },
		 boitef/.style={
        fill=white!70!cyan!70,draw,align=center,
        font=\footnotesize\bfseries,text=black, 
      },
			 boitef coins ronds/.style={boitef,rounded corners=3pt},
			 boitef2/.style={
        fill=white,draw,align=center,
        font=\footnotesize\bfseries,text=white, 
      },
			 boitef2 coins ronds/.style={boitef2,rounded corners=3pt},
			 boitefw/.style={
        fill=white,draw,align=center,
        font=\footnotesize\bfseries,text=black, 
      },
			 boitefw coins ronds/.style={boitefw,rounded corners=3pt},
			 boite2 coins ronds/.style={boite2,rounded corners=3pt},
			boite3 coins ronds/.style={boite3,rounded corners=3pt},
			boite4 coins ronds/.style={boite4,rounded corners=5pt},
      boite coins ronds/.style={boite,rounded corners=3pt},
      boite circulaire/.style={boite4,circle},
      fleche/.style={
        line cap=round,-latex,line width=0.25mm,
        draw=red!60!cyan!60,
      },
			fleche2/.style={
        line cap=round,-latex,line width=0.25mm,
        draw=green!90!cyan!60,
      },
			fleche3/.style={
        line cap=round,-latex,line width=0.25mm,
        draw=blue!90!cyan!60,
      },
			fleche4/.style={
        line cap=round,-latex,line width=0.25mm,
        draw=black!90!cyan!60,
      },
			flechef/.style={
        line cap=round,-latex,line width=0.25mm,
        draw=white!90!cyan!10,
      },
			   grosse flechef/.style={flechef,line width=0.75mm},
      grosse fleche/.style={fleche4,line width=0.75mm},
			grosse fleche2/.style={fleche2,line width=1mm},
			grosse fleche3/.style={fleche3,line width=1mm}
   };

\node[boite4 coins ronds,yshift=17em] (un) {$\underset{x\in S}{\Max}\hspace{0.1cm} f_{1}\left({\color{orange}\overline{x}^{1}},\ldots,\overline{x}^{P}\right)$};

\node[left=1.5cm] (0,0) at (un) {\small\bf Level 1};

\node[boite4 coins ronds,above=0.7cm of un] (un2023){Start};

  \node[boite4 coins ronds,yshift=17em, xshift=12em, text width=2cm]  (a) {Restrict ${\color{orange}\overline{x}^{1}}$ };

\node[boite4 coins ronds,below=of un] (deux) {$\underset{x\in S_{{\bf 1}}}{\Max}\hspace{0.1cm} f_{2}\left({\color{blue}\overline{x}^{1}},{\color{orange}\overline{x}^{2}},\ldots,{\overline{x}^{P}}\right)$};
\node[below=0.32cm] (0,0) at (deux) {\small\bf Level 2};

\node[boite4 coins ronds,left=3cm of deux] (a2) {Restrict \\${\color{orange}\overline{x}^{2}}$ };

\node[boite4 coins ronds,below=of deux] (trois) {$\underset{x\in S_{{\bf 2}}}{\Max}\hspace{0.1cm} f_{3}\left({\color{blue}\overline{x}^{1}},{\color{blue}\overline{x}^{2}},{\color{orange}\overline{x}^{3}},\ldots,{\overline{x}^{P}}\right)$};
\node[right=2cm] (0,0) at (trois) {\small\bf Level 3};

\node[boite4 coins ronds,below=2cm of trois] (quatre) {$\underset{x\in S_{{ P-1}}}{\Max}\hspace{0.1cm} f_{P}\left({\color{blue}\overline{x}^{1}},{\color{blue}\overline{x}^{2}},\ldots,{\color{blue}\overline{x}^{P-1}},{\color{orange}\overline{x}^{P}}\right)$};
\node[right=2.45cm] (0,0) at (quatre) {\small\bf Level P};

	\node[boite2 coins ronds, left=0.2cm of quatre]  (a11f2) {\tiny$S_{{ p-1}}=\left\{x\in S : \text{Restriction on}\hspace{0.13cm}{\color{blue}\overline{x}^{j}} ,j=1,\ldots,{p-1}\right\}$};


\begin{scope}
\draw[grosse fleche] (un2023) -- (un);

\draw[grosse fleche] (un) -- (a);

\draw[grosse fleche] (a) |- (deux);

\draw[grosse fleche] (deux) -- (a2);

\draw[grosse fleche] (a2) |- (trois);

\draw[grosse fleche] (trois) -- (quatre);

		\end{scope}
		
  \end{tikzpicture}
	\caption{Resolution approach for a ML(MO)OLPP}

\end{figure}

\section{Linear Formulation of the $p$th Level}\label{section4}

In this section, we adapt the linear programming method proposed by Sinha and Sinha in \cite{sinha} to express the $p$th level of problem (\ref{mlpp1js}) as a single-objective linear programming problem with bounded variables. This adaptation sets the stage for the subsequent application of the Gabasov's adaptive method algorithm \ref{ama}.

\begin{notation} $\;$\\
\begin{enumerate}
\item  $\stackrel{o}{x}^{p}$ is the optimal solution of $DM_{p}$, and $\stackrel{o}{x}_{i,j}^{p}$ is its $i,j$th component.
\item  $LW_{p-1,j}$ and $UP_{p-1,j}$ are the lower and upper intervals, respectively, of the $(p-1,j)$th component of the decision variable vector.
\end{enumerate}
\end{notation}

\subsection{Construction of the Interval Reduction Map}
The application $\xi_{p-1}$ that we will develop provides the reduced interval of decision variables controlled by $DM_{p-1}$, following the steps (3)-(9) precisely as proposed by Sinha and Sinha.

\begin{notation} $\;$\\
\begin{enumerate}
\item  The superscript in the upper right $^{p}$ indicates that we are in the $p$th level, and the subscript in the lower right $_{i,j}$ designates the $i,j$th component.

\item $l^{(p-1)}$ and $u^{(p-1)}$ represent the lower and upper bounds, respectively, of the decision variables at the $p$th level (these are the new bounds of the decision variables after reducing their intervals by $DM_{p-1}$). Notations $l^{(p-1)}_{i,j}$ and $u^{(p-1)}_{i,j}$ refer to their $i,j$th components.
\end{enumerate}
The index $i=p-1$ signifies that we are addressing decision variables under the control of the $DM_{p-1}$.
\end{notation}

Consider the multilevel linear programming problem (\ref{mlpp1js}), where $P\geq2$. It includes $m$ linear constraints, $n$ decision variables, and a single objective function at each level. For $p=1,\ldots, P$, let $\stackrel{o}{x}^{p}$ be the optimal solution of the following monoobjective linear programming problem:
\begin{equation}\label{Pblm1}
\left\{
\begin{array}{l}
     \underset{x}{\max} \hspace{0.2cm} f_{p}(x)=c_{p}^{T}x           \\
      Ax\leq b\\
x\geq0.
\end{array}
\right.
\end{equation}
Next, we find the bounds for each of the $n$ decision variables as follows:
\begin{equation}\label{boundch5}
l^{(0)}_{i,j}=\underset{q\in\left\{1,\ldots,P\right\}}{\min}\left\{\stackrel{o}{x}_{i,j}^{q}\right\}, \, \, \,
u^{(0)}_{i,j}=\underset{q\in\left\{1,\ldots,P\right\}}{\max}\left\{\stackrel{o}{x}_{i,j}^{q}\right\}, \, \, \, i=1,\ldots, P, j=1,\ldots, n_{i}.
 \end{equation}
Then, set $l^{(1)}_{i,j} \gets l^{(0)}_{i,j}$ and $u^{(1)}_{i,j} \gets u^{(0)}_{i,j}$.  Following this, consider the following $2n$ ideal intervals:
\begin{equation}\label{Pblm2}
\begin{array}{lll}
			{ l^{(1)}_{i,j}}\leq x_{i,j} \leq { u^{(1)}_{i,j}},&\hspace{0.15cm}&i=1,\ldots, P, j=1,\ldots, n_{i}.
			\end{array}
\end{equation}

Guided by $DM_{p}$, each decision variable has the possibility of aligning with an ideal interval or reaching either of its lower or upper bounds. The following are the various scenarios:
\begin{center}
\begin{tikzpicture}
    \tikzset{
      node distance=5mm and 0mm,
      boite4/.style={
        fill=white,draw=black,align=center,
        font=\footnotesize\bfseries,text=black,line width=1pt 
      },
      boite4 coins ronds/.style={boite4,rounded corners=10pt},
      boite4 circulaire/.style={boite4,circle},
			boite4 signal/.style={boite4,signal},
			boite4 circular-sector/.style={boite4,circular sector},
			boite4 chamfered-ellipse/.style={boite4,chamfered ellipse},
			boite4 chamfered-rectangle/.style={boite4,chamfered rectangle},
			boite4 tape/.style={boite4,tape},
			boite4 cylinder/.style={boite4,cylinder},
				boite4 cloud/.style={boite4,cloud},
					boite4 star/.style={boite4,star},
						boite4 ellipse/.style={boite4,ellipse},
    }

    \node[boite4 coins ronds] (deux) {$ l^{(p-1)}_{p-1,j}<\stackrel{o}{x}^{(p-1)}_{p-1,j}<u^{(p-1)}_{p-1,j}$};
    \node[below=0.65cm] at (deux) {\small\bf $1^{\text{st}}$ Case};

    \node[boite4 coins ronds, right=1cm of deux] (a) {$\stackrel{o}{x}^{(p-1)}_{p-1,j}=u^{(p-1)}_{p-1,j}$};
    \node[below=0.65cm] at (a) {\small\bf $2^{\text{nd}}$ Case};

    \node[boite4 coins ronds, left=1cm of deux] (a2) {$\stackrel{o}{x}^{(p-1)}_{p-1,j}=l^{(p-1)}_{p-1,j}$};
    \node[below=0.65cm] at (a2) {\small\bf $3^{\text{rd}}$ Case};

  \end{tikzpicture}
	\end{center}
	
\begin{remark}
Let $\delta_{1}, \delta_{2}, \delta_{3}\in\mathbb{R}$ such that $\delta_{1}<\delta_{3}<\delta_{2}$, and define the applications
\begin{equation}\label{redumap1}
\begin{array}{cccc}
L: & [\delta_{1},\delta_{2}] & \longrightarrow & [\delta_{1},\delta_{3}]\vspace{0.2cm}\\
& x & \longrightarrow & \frac{(\delta_{3}-\delta_{1})x+\delta_{1}(\delta_{2}-\delta_{3})}{\delta_{2}-\delta_{1}}
\end{array},
\end{equation}
\begin{equation}\label{redumapu}
\begin{array}{cccc}
U: & [\delta_{1},\delta_{2}] & \longrightarrow & [\delta_{3},\delta_{2}]\vspace{0.2cm}\\
& x & \longrightarrow & \frac{(\delta_{2}-\delta_{3})x+\delta_{2}(\delta_{3}-\delta_{1})}{\delta_{2}-\delta_{1}}
\end{array}.
\end{equation}
Then $L$ and $U$ are bijective, continuous, and convex mappings that reduce the interval $[\delta_{1},\delta_{2}]$ to $[\delta_{1},\delta_{3}]$ $($resp. $[\delta_{3},\delta_{2}])$. 
This two maps are the key idea behind the construction of the interval reduction map $\xi_{p-1}$.
\end{remark}

\begin{description}[align=left,labelwidth=0.45cm,leftmargin=0cm,labelsep=0cm,itemindent=0cm]
  \item[\textbf{Case 1:}]\hspace{0.05cm} In this case, the $DM_{p-1}$ selects the reduced intervals of the decision variables by examining the sign before the decision variable $x_{p-1,j}$, $j=1,\ldots,n_{p-1}$, in $f_{p-1}(x)$.

If the sign is negative, he  choose the lower interval $LW_{p-1,j}=\left[l_{p-1,j}^{(p-1)},\stackrel{o}{x}^{(p-1)}_{p-1,j}\right]$ using the application:
\begin{equation}
       \begin{array}{llll}
       L_{p-1,j}:& \left[l_{p-1,j}^{(p-1)},u_{p-1,j}^{(p-1)}\right] &  \rightarrow & LW_{p-1,j}            \\
                                                 &  x_{p-1,j}                 &       \rightarrow & \frac{\left(\stackrel{o}{x}^{(p-1)}_{p-1,j}-l_{p-1,j}^{(p-1)}\right)x_{p-1,j}+l_{p-1,j}^{(p-1)}\left(u_{p-1,j}^{(p-1)}-\stackrel{o}{x}^{(p-1)}_{p-1,j}\right)}{u_{p-1,j}^{(p-1)}-l_{p-1,j}^{(p-1)}}
        \end{array}.
\end{equation}
If the sign is positive, he choose the upper interval $UP_{p-1,j} =\left[\stackrel{o}{x}^{(p-1)}_{p-1,j},u_{p-1,j}^{(p-1)}\right]$ using the application:
\begin{equation}
\begin{array}{llll}
U_{p-1,j}:& \left[l_{p-1,j}^{(p-1)},u_{p-1,j}^{(p-1)}\right] & \rightarrow & UP_{p-1,j} \\
                                                 &  x_{p-1,j}                & \rightarrow & \frac{\left(u_{p-1,j}^{(p-1)}-\stackrel{o}{x}^{(p-1)}_{p-1,j}\right)x_{p-1,j}+u_{p-1,j}^{(p-1)}\left(\stackrel{o}{x}^{(p-1)}_{p-1,j}-l_{p-1,j}^{(p-1)}\right)}{u_{p-1,j}^{(p-1)}-l_{p-1,j}^{(p-1)}}
\end{array}.
\end{equation}

Note that applications $L_{p-1,j}$ and $U_{p-1,j}$ are defined similarly to those presented in (\ref{redumap1}) and (\ref{redumapu}), with the following parameters:
$
\begin{array}{ccc}
\delta_{1}=l_{p-1,j}^{(p-1)},& \delta_{2}=u_{p-1,j}^{(p-1)},& \delta_{3}=\stackrel{o}{x}^{(p-1)}_{p-1,j}.
\end{array}
$
These applications specifically affect the $(p-1,j$)th component.

\begin{definition}[Sign function]
The "sign" function assigns the value +1 to positive numbers, -1 to negative numbers, and 0 to zero.
\end{definition}

\begin{proposition}[Definition]\label{def1}For a given decision variable $x_{p-1,j}$, where $x_{p-1,j}\in\left[l_{p-1,j}^{(p-1)},u_{p-1,j}^{(p-1)}\right]$, we define the following quantities:
$$
\begin{array}{ccc}
 T_{p-1,j}^{-}=1-\sign\left(c^{p-1,j}_{p-1}\right)  & \text{and} &T_{p-1,j}^{+}=1+\sign\left(c^{p-1,j}_{p-1}\right).
\end{array}
$$
Then, we define the function $\psi_{p-1,j}\left(x_{p-1,j}\right)$ as follows:
$$
\psi_{p-1,j}\left(x_{p-1,j}\right)=
\frac{T_{p-1,j}^{-}L_{p-1,j}\left(x_{p-1,j}\right)
+T_{p-1,j}^{+}U_{p-1,j}\left(x_{p-1,j}\right)}{2},
$$
when $c^{p-1,j}_{p-1}\neq0$, we get:
$$
\psi_{p-1,j}\left(\left[l_{p-1,j}^{(p-1)},u_{p-1,j}^{(p-1)}\right]\right)=\left\{
\begin{array}{lll}
LW_{p-1,j}&\text{if}&c^{p-1,j}_{p-1}<0,\vspace{0.3cm}\\
UP_{p-1,j}&\text{if}&c^{p-1,j}_{p-1}>0.
\end{array}\right.
$$
\end{proposition}

\begin{proof}\hspace{0.1cm}Since $c^{p-1,j}_{p-1}\neq0$, this means that $\left(c^{p-1,j}_{p-1}>0 \hspace{0.15cm} \text{or}\hspace{0.15cm} c^{p-1,j}_{p-1}<0\right)$. Therefore, we have two cases to consider:
\begin{enumerate}
\item\hspace{0.07cm}  When $c^{p-1,j}_{p-1}>0$, we get:
$
\begin{array}{ccc}
T_{p-1,j}^{+}=2&\text{and}&T_{p-1,j}^{-}=0.
\end{array}
$ 
Then, the function $\psi_{p-1,j}$ becomes:
$$
\psi_{p-1,j}\left(\left[l_{p-1,j}^{(p-1)},u_{p-1,j}^{(p-1)}\right]\right)=U_{p-1,j}\left(\left[l_{p-1,j}^{(p-1)},u_{p-1,j}^{(p-1)}\right]\right)=UP_{p-1,j}.
$$
\item\hspace{0.07cm} When $c^{p-1,j}_{p-1}<0$, we get:
$$
\begin{array}{ccc}
T_{p-1,j}^{+}=0&\text{and}&T_{p-1,j}^{-}=2.
\end{array}
$$
The function $\psi_{p-1,j}$ becomes:
$$
\psi_{p-1,j}\left(\left[l_{p-1,j}^{(p-1)},u_{p-1,j}^{(p-1)}\right]\right)=L_{p-1,j}\left(\left[l_{p-1,j}^{(p-1)},u_{p-1,j}^{(p-1)}\right]\right)=LW_{p-1,j}.
$$
\end{enumerate}
Proposition \ref{c=0} explicitly addresses the situation when $c^{p-1,j}_{p-1}=0$.
\end{proof}

\begin{remark}
We can observe that if $\sign\left(c^{p-1,j}_{p-1}\right)\neq0$, then $\psi_{p-1,j}$ maps the interval $\left[l_{p-1,j}^{(p-1)},u_{p-1,j}^{(p-1)}\right]$ either to $LW_{p-1,j}$ or to $UP_{p-1,j}$. However, if $\sign\left(c^{p-1,j}_{p-1}\right)=0$, it means that the component $x_{p-1,j}$ does not appear in the function $f_{p-1}$ and thus does not influence the maximization of this function. In this case, we do not reduce the interval $\left[l_{p-1,j}^{(p-1)},u_{p-1,j}^{(p-1)}\right]$ because it would have no effect on the function's maximization. We simply apply the identity function to keep the interval unchanged.

It is important to note that these updates will only be applied if $c^{p-1,j}_{p}\neq0$. Otherwise, the $DM_{p}$ at the lower levels will not use them because the relevant variable does not appear in the function $f_{p}$.

All these considerations are taken into account by the application $\nu_{p-1,j}$, which handles the case where $c^{p-1,j}_{p}=0$, as indicated in proposition \ref{c=0}.
\end{remark}

\begin{proposition}[Definition]\label{c=0}For a given decision variable $x_{p-1,j}$, where $x_{p-1,j}\in\left[l_{p-1,j}^{(p-1)},u_{p-1,j}^{(p-1)}\right]$, we define the following quantities:
$$
\begin{array}{ccc}
 S_{p-1}^{p-1,j}=\sign\left(c^{p-1,j}_{p-1}\right) &
 \text{and}&
 S_{p}^{p-1,j}=\sign\left(c^{p-1,j}_{p}\right).
\end{array}
$$
Using the notation from Proposition \ref{c=0}, we then define:
$$
\nu_{p-1,j}\left(x_{p-1,j}\right)=\left(S_{p}^{p-1,j}S_{p-1}^{p-1,j}\right)^{2}\psi_{p-1,j}\left(x_{p-1,j}\right)
+\left(1-\left(S_{p}^{p-1,j}S_{p-1}^{p-1,j}\right)^{2}\right)x_{p-1,j}.
$$
Thus, if both $c^{p-1,j}_{p-1}\neq 0$ and $c^{p-1,j}_{p}\neq0$, we have:
$$
\nu_{p-1,j}\left(\left[l_{p-1,j}^{(p-1)},u_{p-1,j}^{(p-1)}\right]\right)=\psi_{p-1,j}\left(\left[l_{p-1,j}^{(p-1)},u_{p-1,j}^{(p-1)}\right]\right).
$$
Else, if $c^{p-1,j}_{p-1}=0$ or $c^{p-1,j}_{p}=0$, then we simply apply the identity application, that is:
$
\nu_{p-1,j}=Id_{\left[l_{p-1,j}^{(p-1)},u_{p-1,j}^{(p-1)}\right]}.
$
\end{proposition}

\begin{proof}$\;$\\
\begin{enumerate}
\item\hspace{0.07cm} When both $c^{p-1,j}_{p-1}\neq 0$ and $c^{p-1,j}_{p}\neq0$, we have:
$$
\begin{array}{ccc}
(S_{p}^{p-1,j}S_{p-1}^{p-1,j})^{2}=1, &1-(S_{p}^{p-1,j}S_{p-1}^{p-1,j})^{2}=0
\end{array}
$$
and
$$
\nu_{p-1,j}\left(\left[l_{p-1,j}^{(p-1)},u_{p-1,j}^{(p-1)}\right]\right)=\psi_{p-1,j}\left(\left[l_{p-1,j}^{(p-1)},u_{p-1,j}^{(p-1)}\right]\right).
$$
\item\hspace{0.07cm} Otherwise, if $c^{p-1,j}_{p-1}=0$ or $c^{p-1,j}_{p}=0$, then we simply apply the identity function, which means:
$$
\begin{array}{ccc}
(S_{p}^{p-1,j}S_{p-1}^{p-1,j})^{2}=0&\text{and}&1-(S_{p}^{p-1,j}S_{p-1}^{p-1,j})^{2}=1
\end{array}
$$
and
$$
\nu_{p-1,j}\left(\left[l_{p-1,j}^{(p-1)},u_{p-1,j}^{(p-1)}\right]\right)=\left[l_{p-1,j}^{(p-1)},u_{p-1,j}^{(p-1)}\right]=Id_{\left[l_{p-1,j}^{(p-1)},u_{p-1,j}^{(p-1)}\right]}.
$$
\end{enumerate}
\end{proof}

  \item[\textbf{Case 2:}] \hspace{0.05cm} If $\stackrel{o}{x}^{(p-1)}_{p-1,j}=u^{(p-1)}_{p-1,j}$, then the choice of interval reduction will be as follows:
$$
LW_{p-1,j}^{{\color{red}\alpha_{p-1,j}}}=\left[l_{p-1,j}^{(p-1)},u_{p-1,j}^{(p-1)}-{\color{red}\alpha_{p-1,j}}\right],
$$
using the application:
\begin{equation}
\begin{array}{llll}
L_{p-1,j}^{{\color{red}\alpha_{p-1,j}}}:& \left[l_{p-1,j}^{(p-1)},u_{p-1,j}^{(p-1)}\right] & \rightarrow & LW_{p-1,j}^{{\color{red}\alpha_{p-1,j}}}\vspace{0.15cm}\\
                                                 &  x_{p-1,j}                         & \rightarrow & \frac{\left(u_{p-1,j}^{(p-1)}-\left(l_{p-1,j}^{(p-1)}+{\color{red}\alpha_{p-1,j}}\right)\right)x_{p-1,j}+{\color{red}\alpha_{p-1,j}} l_{p-1,j}^{(p-1)}}{u_{p-1,j}^{(p-1)}-l_{p-1,j}^{(p-1)}}
\end{array}.
\end{equation}
where ${\color{red}\alpha_{p-1,j}}\in\mathbb{R}$ is an arbitrary number chosen by $DM_{p-1}$, such that:
\begin{equation}\label{alpha}
{\color{red}\alpha_{p-1,j}}<u_{p-1,j}^{(p-1)}-l_{p-1,j}^{(p-1)}.
\end{equation}
\item[\textbf{Case 3:}]  \hspace{0.05cm} If $\stackrel{o}{x}^{p-1,j}_{p-1}=l^{(p-1)}_{p-1,j}$, then the choice of interval reduction will be as follows:
$$
UP_{p-1,j}^{{\color{red}\alpha_{p-1,j}}}=\left[l_{p-1,j}^{(p-1)}+{\color{red}\alpha_{p-1,j}},u_{p-1,j}^{(p-1)}\right],
$$
using the application:
\begin{equation}
\begin{array}{llll}
U_{p-1,j}^{{\color{red}\alpha_{p-1,j}}}:& \left[l_{p-1,j}^{(p-1)},u_{p-1,j}^{(p-1)}\right] & \rightarrow & UP_{p-1,j}^{{\color{red}\alpha_{p-1,j}}}\vspace{0.15cm}\\
                                                 &  x_{p-1,j}                         & \rightarrow & \frac{\left(u_{p-1,j}^{(p-1)}-\left(l_{p-1,j}^{(p-1)}+{\color{red}\alpha_{p-1,j}}\right)\right)x_{p-1,j}+{\color{red}\alpha_{p-1,j}}u_{p-1,j}^{(p-1)}}{u_{p-1,j}^{(p-1)}-l_{p-1,j}^{(p-1)}}
\end{array}.
\end{equation}


\begin{proposition}[Definition]\label{def2}
For each component $\stackrel{o}{x}^{(p-1)}_{p-1,j}$, we define the following quantities:
$$
\begin{array}{ccc}
A_{p-1,j}=\sign\left(u^{(p-1)}_{p-1,j}-\stackrel{o}{x}^{(p-1)}_{p-1,j}\right)& \text{and} &B_{p-1,j}=\sign\left(\stackrel{o}{x}^{(p-1)}_{p-1,j}-l^{(p-1)}_{p-1,j}\right).
\end{array}
$$
Then, we define the function $\hat{\psi}_{p-1,j}$ as follows:
$$
\hat{\psi}_{p-1,j}(x_{p-1,j})=B_{p-1,j}L_{p-1,j}^{{\color{red}\alpha_{p-1,j}}}(x_{p-1,j})
+A_{p-1,j}U_{p-1,j}^{{\color{red}\alpha_{p-1,j}}}(x_{p-1,j}).
$$ 
For any $x_{p-1,j}\in\left[l_{p-1,j}^{(p-1)},u_{p-1,j}^{(p-1)}\right]$, the following cases apply:
\begin{itemize}
\item[$\diamond$]\hspace{0.07cm}  If $A_{p-1,j}=0$, then $\hat{\psi}_{p-1,j}\left(x_{p-1,j}\right)=L_{p-1,j}^{{\color{red}\alpha_{p-1,j}}}(x_{p-1,j})$.
\item[$\diamond$]\hspace{0.07cm}   If $B_{p-1,j}=0$, then $\hat{\psi}_{p-1,j}\left(x_{p-1,j}\right)=U_{p-1,j}^{{\color{red}\alpha_{p-1,j}}}(x_{p-1,j})$.
\end{itemize}
\end{proposition}

\begin{proof}\hspace{0.07cm} 
This proposition is self-evident.
\end{proof}

\begin{proposition}[Definition]\label{c=1}Let $x_{p-1,j}\in\left[l_{p-1,j}^{(p-1)},{u_{p-1,j}^{(p-1)}}\right]$, using all the previous notations, we define:
$$
\hat{\nu}_{p-1,j}\left(x_{p-1,j}\right)=\left(S_{p}^{p-1,j}S_{p-1}^{p-1,j}\right)^{2}\hat{\psi}_{p-1,j}\left(x_{p-1,j}\right)+\left(1-\left(S_{p}^{p-1,j}S_{p-1}^{p-1,j}\right)^{2}\right)x_{p-1,j}.
$$
Then,
$$
\hat{\nu}_{p-1,j}\left(\left[l_{p-1,j}^{(p-1)},{u_{p-1,j}^{(p-1)}}\right]\right)=\left\{\begin{array}{ll}
\hat{\psi}_{p-1,j}\left(\left[l_{p-1,j}^{(p-1)},{u_{p-1,j}^{(p-1)}}\right]\right)& \text{if}\hspace{0.12cm}c^{p-1,j}_{p-1}\neq 0 \hspace{0.15cm}\text{and}\hspace{0.15cm} c^{p-1,j}_{p}\neq0,\vspace{0.17cm}\\
\hat{\nu}_{p-1,j}=Id_{\left[l_{p-1,j}^{(p-1)},{u_{p-1,j}^{(p-1)}}\right]}&\text{else}.
\end{array}\right.
$$
\end{proposition}

\begin{proof}\hspace{0.07cm}
The proof is similar to that of Proposition \ref{c=0}.
\end{proof}
\end{description}


\begin{notation}
Using all the previous notations, we denote the following:
$$
\begin{array}{llll}
1.\hspace{0.2cm}A^{+}_{p-1,j}=1+A_{p-1,j},&2.\hspace{0.2cm}A^{-}_{p-1,j}=1-A_{p-1,j},&3.\hspace{0.2cm}B^{-}_{p-1,j}=1-B_{p-1,j},& \vspace{0.2cm}\\
4.\hspace{0.2cm}B^{+}_{p-1,j}=1+B_{p-1,j},&5.\hspace{0.2cm}H_{p-1,j}^{1}=A^{+}_{p-1,j}A^{-}_{p-1,j}+B^{+}_{p-1,j}B^{-}_{p-1,j}, &6.\hspace{0.2cm}H_{p-1,j}^{2}=-A_{p-1,j}B_{p-1,j}.&
\end{array}
$$
\end{notation}

	\begin{remark}
In practice, in the 1$st$ iteration at the 1$st$ level, $DM_{1}$ solves its problem with respect to the initial constraints only, without any reduced interval of its variables. After obtaining its solution, it performs the reduction, which will be practically carried out in the second iteration by $DM_2$. Then, we define the interval reduction application at this level as follows: $\xi_{0}\left(x_{i,j}\right)=x_{i,j}$.
	\end{remark}
\begin{definition}[Interval Reduction Map (IRM)]\label{RRM}
Let $i = 2, \ldots, P$, $j = 1, \ldots, n_{i}$, and $x_{i,j} \in \left[l^{(p-1)}{i,j}, u^{(p-1)}{i,j}\right]$. we define:
 \begin{equation}
\xi_{p-1}(x_{i,j})=   \left\{\begin{array}{lll}
                              H_{p-1,j}^{1}\hat{\psi}_{i,j}(x_{i,j})+H_{p-1,j}^{2}\hat{\nu}_{i,j}(x_{i,j})&\text{si}&i=p-1 \\
                               x_{i,j}                                                                                    &\text{si} & i\neq p-1.
                        \end{array}
                        \right.
 \end{equation}
\end{definition}

\begin{proposition}
The $\xi_{p-1}$ application will adjust the decision variable intervals $x_{i,j}$ by taking into account all the modifications made previously.
\end{proposition}

\begin{proof}Let  $i=2,\ldots,P$, $j=1,\ldots,n_{i}$ and $ x_{i,j}\in\left[l^{(p-1)}_{i,j},u^{(p-1)}_{i,j}\right]$, then
\begin{description}[align=left,labelwidth=0cm,leftmargin=0cm,labelsep=0cm,itemindent=0cm]
\item[Proof of the first case:] \hspace{0.07cm}When $l^{(p-1)}_{p-1,j}<\stackrel{o}{x}^{(p-1)}_{p-1,j}<u^{(p-1)}_{p-1,j}$, we get:
$$
\begin{array}{c}
\shadowbox{$\begin{array}{l}
A_{p-1,j}=1,\\
B_{p-1,j}=-1
\end{array}$}\Rightarrow\left\{
\begin{array}{l}
\begin{array}{lcl}
 A^{+}_{p-1,j}=1+A_{p-1,j}=2,&&A^{-}_{p-1,j}=1-A_{p-1,j}=0,\vspace{0.2cm}\\
B^{+}_{p-1,j}=1+B_{p-1,j}=0,&&B^{-}_{p-1,j}=1-B_{p-1,j}=2,
\end{array}\vspace{0.2cm}\\
\begin{array}{lclcl}
H_{p-1,j}^{1}=A^{+}_{p-1,j}A^{-}_{p-1,j}+B^{+}_{p-1,j}B^{-}_{p-1,j}&=&2\times 0+0\times 2&=&0, \vspace{0.2cm}\\
 H_{p-1,j}^{2}=-A_{p-1,j}B_{p-1,j}&=&-1\times -1&=&1.
\end{array}
\end{array}\right.
\end{array}
$$
Then, we obtain:
$$
\xi_{p-1}(x_{i,j})=   \left\{\begin{array}{lclcl}
                              0\times\hat{\psi}_{i,j}(x_{i,j})+\hat{\nu}_{i,j}(x_{i,j})&=&\hat{\nu}_{i,j}(x_{i,j})&\text{if}&i=p-1, \vspace{0.2cm}\\
                               x_{i,j}         &&                                                                           &\text{if} & i\neq p-1.
                        \end{array}
                        \right.
$$

  \item[Proof of the second case:] \hspace{0.07cm}When $\stackrel{o}{x}^{(p-1)}_{p-1,j}=u^{(p-1)}_{p-1,j}$, we get:

$$
\begin{array}{c}
\shadowbox{$\begin{array}{l}
A_{p-1,j}=0,\\
B_{p-1,j}=1
\end{array}$}\Rightarrow\left\{\begin{array}{l}
\begin{array}{lcl}
 A^{+}_{p-1,j}=1+A_{p-1,j}=1,&&A^{-}_{p-1,j}=1-A_{p-1,j}=1,\vspace{0.2cm}\\
B^{+}_{p-1,j}=1+B_{p-1,j}=2,&&B^{-}_{p-1,j}=1-B_{p-1,j}=0,
\end{array}\vspace{0.2cm}\\
\begin{array}{lclcl}
H_{p-1,j}^{1}=A^{+}_{p-1,j}A^{-}_{p-1,j}+B^{+}_{p-1,j}B^{-}_{p-1,j}&=&1\times 1+2\times 0&=&1, \vspace{0.2cm}\\
 H_{p-1,j}^{2}=-A_{p-1,j}B_{p-1,j}&=&-0\times 1&=&0.
\end{array}
\end{array}\right.
\end{array}
$$
Then, we obtain:
$$
\xi_{p-1}(x_{i,j})=   \left\{\begin{array}{lclcl}
                              \hat{\psi}_{i,j}(x_{i,j})+0\times\hat{\nu}_{i,j}(x_{i,j})&=&\hat{\psi}_{i,j}(x_{i,j})&\text{if}&i=p-1, \vspace{0.2cm}\\
                               x_{i,j}         &&                                                                           &\text{if} & i\neq p-1.
                        \end{array}
                        \right.
$$

  \item[Proof of the third case:] \hspace{0.07cm}When $\stackrel{o}{x}^{(p-1)}_{p-1,j}=l^{(p-1)}_{p-1,j}$, we get:
	
$$
\begin{array}{c}
\shadowbox{$\begin{array}{l}
A_{p-1,j}=1,\\
B_{p-1,j}=0
\end{array}$}\Rightarrow\left\{
\begin{array}{l}
\begin{array}{lcl}
A^{+}_{p-1,j}=1+A_{p-1,j}=2,&&A^{-}_{p-1,j}=1-A_{p-1,j}=0,\vspace{0.2cm}\\
B^{+}_{p-1,j}=1+B_{p-1,j}=1,&&B^{-}_{p-1,j}=1-B_{p-1,j}=1,
\end{array}\vspace{0.2cm}\\
\begin{array}{lclcl}
H_{p-1,j}^{1}=A^{+}_{p-1,j}A^{-}_{p-1,j}+B^{+}_{p-1,j}B^{-}_{p-1,j}&=&2\times 0+1\times 1&=&1, \vspace{0.2cm}\\
 H_{p-1,j}^{2}=-A_{p-1,j}B_{p-1,j}&=&-1\times 0&=&0.
\end{array}
\end{array}\right.
\end{array}
$$
Then, we obtain:
$$
\xi_{p-1}(x_{i,j})=   \left\{\begin{array}{lclcl}
                              \hat{\psi}_{i,j}(x_{i,j})+0\times\hat{\nu}_{i,j}(x_{i,j})&=&\hat{\psi}_{i,j}(x_{i,j})&\text{if}&i=p-1, \vspace{0.2cm}\\
                               x_{i,j}         &&                                                                           &\text{if} & i\neq p-1.
                        \end{array}
                        \right.
$$
\end{description}\end{proof}

\begin{remark}
The application $\xi_{p-1}$ will reduce the intervals of the components $x_{p-1,j}$ for all $j=1,\ldots,n_{p}$ only.
\end{remark}

\begin{notation}
We denote by $l^{(p)}$ and $u^{(p)}$ the following vectors:
\begin{equation}\label{equation1}
l^{(p)}:=\xi_{p-1}(l^{(p-1)})=\left(\xi_{p-1}\left(l^{(p-1)}_{i,j}\right), i=1,\ldots,P, j=1,\ldots,n_{i}\right),
\end{equation}
\begin{equation}\label{equation1}
u^{(p)}:=\xi_{p-1}(u^{(p-1)})=\left(\xi_{p-1}\left(u^{(p-1)}_{i,j}\right), i=1,\ldots,P, j=1,\ldots,n_{i}\right).
\end{equation}
\end{notation}

\begin{definition}
The linear formulation of the $p$th level is given as follows:
	\begin{equation}\label{model}
	\left\{
\begin{array}{l}
	            \underset{x}{\max}\hspace{0.12cm} f_{p}(x)=c^{T}_{p}x\\
                                             Ax\leq b                 \\
		                                         \xi_{p-1}\left(l^{(p-1)}\right)\leq x\leq \xi_{p-1}\left(u^{(p-1)}\right).
		                        
\end{array}\right.
\end{equation}
\end{definition}

	\section{Algorithm to solve the ML(MO)OLPP (\ref{mlpp1js}) using Gabasov's adaptive method algorithm \ref{ama}}\label{section5}
The algorithm follows a systematic procedure. Initially, individual solutions for each objective function are computed independently, providing insight into the optimization landscape for each criterion. Following this, careful consideration is given to establishing bounds for all decision variables, a crucial step that sets the stage for subsequent computations. The compromise search phase is then initiated, involving the solution of $P-1$ standard linear programming problems with bounded variables using the adaptive method algorithm \ref{ama}.
\par Consider the following linear programming problem:
	\begin{equation}\label{model2}
	\left\{
\begin{array}{l}
	            \underset{x}{\max}\hspace{0.15cm} f_{p}(x)=c^{T}_{p}x\\
                                             Ax\leq b                 \\
		                                         l^{(p)}\leq x\leq u^{(p)}.
\end{array}\right.
\end{equation}
\begin{algorithm}[H]
\caption{Procedure for resolving the ML(MO)OLPP\label{alg:mlmolpp}}

\textcolor{blue}{\textbf{Step 1:}} Set $p=1$.\\
{\bf While $p \leq P$, repeat the following steps:}
\begin{enumerate}
\item  Explicit the model (\ref{Pblm1}).
\item  Obtain the solution $\stackrel{o}{x}^{p}$ of the model (\ref{Pblm1}) "optimal solution of the $DM_{p}$".
\item  Set $p \gets p+1$.
\end{enumerate}
{\bf End while}\\
\textcolor{blue}{\textbf{Step 2:}} Obtain all bounds $l^{(1)}_{i,j}$ and $u^{(1)}_{i,j}$ using formula (\ref{boundch5}).\\
\textcolor{blue}{\textbf{Step 3:}} Set $p \gets 2$.\\
{\bf While $p \leq P$, repeat the following steps:}
\begin{enumerate}
\item  Explicit the model (\ref{model2}).
\item  Solve the model (\ref{model2}) using the Gabasov's monocriteria adaptive method algorithm \ref{ama} to obtain the optimal solution $\stackrel{c}{x}^{p}$.
\item  Set $p \gets p+1$.
\end{enumerate}
{\bf End while}\\
\textcolor{blue}{\textbf{Output:}} The satisfactory solution $\stackrel{c}{x}^{P}$ for the multilevel problem.
\end{algorithm}

\section{Numerical Example: Vaccination Planning in Long-Term Care Facilities}\label{section6}
Vaccination planning in long-term care facilities (LTCFs) presents a formidable challenge, demanding meticulous consideration of intricate variables such as facility capacity, demographics of the target population, and the specific requirements of vulnerable cohorts. This study is focused on elucidating the complexities inherent in the United Kingdom (UK) LTCF vaccination strategy for the year 2021. A mathematical model is deployed to simulate the distribution of COVID-19 vaccines across varied facilities, systematically incorporating constraints such as capacity limitations and the presence of vulnerable populations. The data employed in the model are derived from publicly available sources.

\begin{sidewaystable}[!htbp]

  \centering
\begin{tabular}{
    l
    S[table-format=2.2]
    S[table-format=1.2]
    S[table-format=1.2]
    S[table-format=1.2]
    S[table-format=4.0]
    S[table-format=1.1]
    S[table-format=1.1]
    S[table-format=1.1] 
  }
    \toprule
   {\textbf{Region} ($i$) } & {\textbf{Population} ($P_{i}$)}  & {\textbf{Cases} ($I_{i}$)}  & {\textbf{Administration doses Capacity} ($C_{i}$)}  & {\textbf{Equitable Vaccine Distribution} ($E_{i}$)}\\
    {\centering} & { (millions)} & { (millions) } & {(millions/year)} & {(\%)} \\
    \midrule
    \textit{1-England} & 56.48 & 8.55 &51.2& 80 \\
    \textit{2-Scotland} & 5.45 & 0.45 &4.95& 70 \\
    \textit{3-Wales} & 3.16 & 0.15 &2.87& 60 \\
    \textit{4-Northern Ireland} & 1.91 & 0.1 &1.78& 50 \\
    \bottomrule
  \end{tabular}
  \vspace{\baselineskip}

  \begin{tabular}{llcc}
    \toprule
    \textbf{Region} ($i$) & \textbf{Hospital} ($j$) & \textbf{Vulnerable Population} ($V_{ij}$) & \textbf{Care Hospital Capacity} ($B_{ij}$) \\
    & & (millions) & beds (millions)  \\
    \midrule
    {\it1-England} & {\it1-Royal Free Hospital (RFH)} & 0.025 & 0.030 \\
    & {\it2-University College Hospital (UCLH)} & 0.030 & 0.035 \\
    & {\it3-University Hospitals of Leicester (UHL)} & 0.020 & 0.025 \\
    & {\it4-Addenbrooke's Hospital (ADH)} & 0.028 & 0.032 \\
    & {\it5-Guy's and St Thomas' NHS Foundation Trust (GSTT)} & 0.032 & 0.034 \\
    & {\it6-Imperial College Healthcare NHS Trust (ICH)} & 0.027 & 0.030 \\
    & {\it7-University Hospitals Birmingham NHS Trust (UHB)} & 0.035 & 0.038 \\
    & {\it8-Manchester University NHS Foundation Trust (MFT)} & 0.023 & 0.028 \\
    \midrule
    {\it2-Scotland} & {\it1-Queen Elizabeth University Hospital (QEH)} & 0.035 & 0.040 \\
    & {\it2-Ninewells Hospital and Medical School (9W)} & 0.020 & 0.025 \\
    & {\it3-Glasgow Royal Infirmary (GRI)} & 0.030 & 0.033 \\
    & {\it4-Aberdeen Royal Infirmary (ARI)} & 0.025 & 0.028 \\
    & {\it5-Royal Edinburgh Hospital (REH)} & 0.028 & 0.031 \\
    \midrule
    {\it3-Wales} & {\it1-Ysbyty Gwynedd (YG)} & 0.015 & 0.020 \\
    & {\it2-University Hospital of Wales (UHW)} & 0.012 & 0.015 \\
    & {\it3-Morriston Hospital (MH)} & 0.010 & 0.012 \\
    & {\it4-Royal Gwent Hospital (RGH)} & 0.008 & 0.010 \\
    & {\it5-Wrexham Maelor Hospital (WMH)} & 0.007 & 0.008 \\
    \midrule
    {\it4-Northern Ireland} & {\it1-Belfast City Hospital (BCH)} & 0.018 & 0.019 \\
    & {\it2-Altnagelvin Area Hospital (AAH)} & 0.016 & 0.017 \\
    & {\it3-Craigavon Area Hospital (CAH)} & 0.014 & 0.015 \\
    & {\it4-South West Acute Hospital (SWAH)} & 0.012 & 0.013 \\
    \bottomrule
  \end{tabular}
\vspace{0.1cm}
  \caption{Comprehensive Real-World Data for UK: Population, COVID-19 Cases, Vaccination Capacity, and Healthcare Infrastructure in 2021\label{tab:real-world-data}}
\end{sidewaystable}

\subsection{Model Formulation \label{MDLF}}

The analysis is regionally structured, with the UK segregated into four distinct administrative regions: England (Region 1), Scotland (Region 2), Wales (Region 3), and Northern Ireland (Region 4). Each region, denoted by $i = 1, 2, 3, 4$, is characterized by unique demographic and healthcare parameters. A detailed examination of hospitals within each region is undertaken, integrating information on vulnerable populations and care unit capacities, as presented in Table \ref{tab:real-world-data}.

For clarity, in Region 1 (England), hospitals are labeled as $11$ (Royal Free Hospital - RFH), $12$ (University College Hospital - UCLH), $13$ (University Hospitals of Leicester - UHL), $14$ (Addenbrooke's Hospital - ADH),  $15$ (Guy's and St Thomas' NHS Foundation Trust - GSTT), $16$ (Imperial College Healthcare NHS Trust - ICH), $17$ (University Hospitals Birmingham NHS Trust - UHB), $18$ (Manchester University NHS Foundation Trust - MFT), and so on, with $j_{1}=8$.

In Region 2 (Scotland), hospitals are denoted as $21$ (Queen Elizabeth University Hospital - QEH), $22$ (Ninewells Hospital and Medical School - 9W), $23$ (Glasgow Royal Infirmary - GRI), $24$ (Aberdeen Royal Infirmary - ARI), $25$ (Royal Edinburgh Hospital - REH), and so on, with $j_{2}=5$.

Region 3 (Wales) is represented by hospital $31$ (Ysbyty Gwynedd - YG), $32$ (University Hospital of Wales - UHW), $33$ (Morriston Hospital - MH), $34$ (Royal Gwent Hospital - RGH), $35$ (Wrexham Maelor Hospital - WMH), and so on, with $j_{3}=5$.

In Region 4 (Northern Ireland), features hospital $41$  (Belfast City Hospital - BCH), $42$ (Altnagelvin Area Hospital - AAH), $43$ (Craigavon Area Hospital - CAH), $44$ (South West Acute Hospital - SWAH), and so on, with $j_{4}=4$.

Our primary goal is to devise a robust multilevel optimization framework that accurately emulates the intricacies of the COVID-19 vaccine allocation process, systematically addressing the specific priorities at each administrative tier. We meticulously elaborate on the hierarchical structure of the model, elucidating the decision variables at each level, formulating the corresponding objective functions, and systematically enumerating the constraints essential for a faithful representation of the vaccination planning landscape in the UK. Our formulated model aims to provide not only valuable insights but also strategic directives aimed at optimizing the distribution of COVID-19 vaccines within the intricate framework of LTCFs in the UK.

We present a multilevel optimization problem designed to simulate the vaccination allocation process, strategically organized across three distinct levels:
\begin{description}
\item[Level 1: Central Management (Government of the UK)]$\;$\\
\begin{description}
\item[$\diamond$ \bf Objective:] Maximize the total number of COVID-19 vaccine doses allocated to the UK.
\item[$\diamond$ \bf Possible Decisions:] $x_{11}$ - Number of vaccine doses allocated to the UK.
\item[$\diamond$ \bf Objective Function (Level 1):] Maximize $x_{11}$ (vaccine doses allocated to the UK).
\item[$\diamond$ \bf Constraints (Level 1):] National vaccine production capacity: $x_{11} \leq 100$ million doses (assuming a procurement capacity of 100 million doses).
\end{description}

\item[\textbf{Level 2: Regional Management (Regions of the UK)}]$\;$\\
\begin{description}
\item[$\diamond$ \bf Objective:] Optimize vaccine distribution at the regional level based on population, infection rates, and healthcare infrastructure.
\item[$\diamond$ \bf Possible Decisions:] $x_{2i}$ represents the allocation of COVID-19 vaccine doses to National Health Service (NHS) region $i$.
\item[$\diamond$ \bf Objective Function (Level 2):] Maximize $\sum_{i=1}^{4} \lambda_{i}x_{2i}$ (vaccine utilization by region), where $\lambda_{i}$ are positive coefficients with $\sum_{i=1}^{4} \lambda_{i} = 1$.\\

The objective at Level 2 is to optimize the distribution of COVID-19 vaccine doses across the four NHS regions, represented by the objective function $(x_{21}, x_{22}, x_{23}, x_{24})$. The coefficients $\lambda_{i}$ reflect the priority assigned to each region, and they are determined based on the urgency and specific considerations set by the government.
\item[$\diamond$ \bf Constraints (Level 2):]$\;$\\
\begin{itemize}
\item[$\bullet$] Regional population:   $ x_{2i} \leq P_{i}$ (total population of NHS region $i$).
\item[$\bullet$] Regional infection rates: $ x_{2i} \geq I_{i}$ (total number of confirmed COVID-19 cases in NHS region $i$).
\item[$\bullet$] Regional healthcare infrastructure:  $ x_{2i} \geq C_{i}$ (total vaccine administration capacity of NHS region $i$).
\item[$\bullet$] Vaccination allocation to regions constraint: $\sum_{i=1}^{4} x_{2i}\leq x_{11}$.
\end{itemize}
\end{description}

\item[\textbf{Level 3: Local Management (Hospitals within Each Region)}]$\;$\\
\begin{description}
\item[$\diamond$ \bf Objective: ]
Optimize vaccine allocation within each NHS region by prioritizing vulnerable populations and ensuring equitable access to vaccines.
\item[$\diamond$ \bf Possible Decisions: ]
$x_{3i}^{j}$ represents the number of COVID-19 vaccine doses allocated to NHS hospital $j$ in NHS region $i$, where $j=1,\ldots,j_{i}$.
\item[$\diamond$ \bf Objective Function (Level 3): ]
Maximize $\sum_{i=1}^{4}\sum_{j=1}^{j_{i}} \omega_{i}^{j}x_{3i}^{j}$, where $\omega_{i}^{j}=\lambda_{i}\xi_{i}^{j}$, subject to $\sum_{j=1}^{j_{i}}\xi_{i}^{j}=1$ and $\xi_{i}^{j}$ being positive coefficients.\\
The objective at Level 3 is to optimize the allocation of COVID-19 vaccine doses among NHS hospitals within each NHS hospital $j$ in NHS region $i$, represented by the objective function: $$(x_{31}^{1}, x_{31}^{2}, x_{31}^{3},x_{31}^{4},x_{31}^{5},x_{31}^{6},x_{31}^{7},x_{31}^{8},x_{32}^{1},x_{32}^{2},x_{32}^{3},x_{32}^{4},x_{32}^{5},x_{33}^{1},x_{33}^{2},x_{33}^{3},x_{33}^{4},x_{33}^{5},x_{34}^{1},x_{34}^{2},x_{34}^{3},x_{34}^{4}).$$ The coefficients $\omega_{i}^{j}$ reflect the priority assigned to each NHS hospital $j$ in NHS region $i$, and they are determined based on the urgency and specific considerations set by the government. 
\item[$\diamond$ \bf Constraints (Level 3): ]$\;$\\
\begin{itemize}
\item[$\bullet$] Care hospital capacity:  $x_{3i}^{j} \leq  B_{ij}$ (total number of beds in NHS hospital $j$ in NHS region $i$).\\
This constraint ensures that the number of vaccine doses allocated to a hospital does not exceed its capacity to administer them based on the number of beds. This ensures efficient use of resources.
\item[$\bullet$] Vulnerable populations:  $  x_{3i}^{j} \geq V_{ij}$  (number of individuals in vulnerable groups at NHS hospital $j$ in NHS region $i$).\\
This constraint ensures that the number of vaccine doses allocated to a hospital is sufficient to cover at least the number of vulnerable individuals within its patient population. This prioritizes vulnerable groups.
\item[$\bullet$] Vaccination allocation to hospitals constraint: $\sum_{i=1}^{4}\sum_{j=1}^{j_{i}} x_{3i}^{j} \leq \sum_{i=1}^{4} x_{2i}$ (total number of vaccine doses allocated to the four regions).\\
This constraint ensures that the total number of vaccine doses allocated to hospitals in all regions does not exceed the total allocated to the UK. This maintains consistency and prevents exceeding the available supply.
\end{itemize}
\end{description}
\end{description}

\begin{table}[htb]
		\centering
	\renewcommand{\arraystretch}{1} 
  
  \begin{tabular}{c|c|l}
    \toprule
    \textbf{Levels}  & \textbf{Objective functions} & \textbf{Constraints} \\
    \midrule
    \textit{Level 1} & $f_1(x) = x_{11}$ & $x_{11} \leq 100$ \\
    \midrule
		\textit{Level 2} &  $f_2(x) = \sum_{i=1}^{4} \lambda_{i}x_{2i}$ &
		\begin{tabular}{l}
  $
    \left\{
      \begin{array}{lcl}
        \diamond \quad  x_{2i} \leq P_{i}&\Leftrightarrow & \mathbf{Id}_{4}(\overline{x}^{2})^{T}\leq  \mathbf{P}\vspace{0.25cm} \\
        \diamond \quad  x_{2i} \geq I_{i}&\Leftrightarrow &-\mathbf{Id}_{4}(\overline{x}^{2})^{T}\leq -\mathbf{I}\vspace{0.25cm} \\
        \diamond \quad  x_{2i} \geq C_{i}&\Leftrightarrow &-\mathbf{Id}_{4}(\overline{x}^{2})^{T}\leq -\mathbf{C}\vspace{0.25cm} \\
        \diamond \quad \sum_{i=1}^{4} x_{2i} \leq x_{11}&\Leftrightarrow &\mathbf{1}_{4}(-\overline{x}^{1},\overline{x}^{2})^{T}\leq 0
      \end{array}
    \right.
  $\\
	\\
 \vspace{0.5cm}
  \fbox{
  
      \textit{13 Constraints}
  }
\end{tabular}
	
\\
    \midrule
    \textit{Level 3} & $f_3(x) = \sum_{i=1}^{4}\sum_{j=1}^{j_{i}} \omega_{i}^{j}x_{3i}^{j}$ &
		\begin{tabular}{l}
		$
\left\{
  \begin{array}{lcl}
    \diamond  \quad x_{3i}^{j} \leq B_{ij}&\Leftrightarrow &\mathbf{Id}_{7}(\overline{x}^{3})^{T}\leq \mathbf{B}\vspace{0.25cm}\\
    \diamond \quad x_{3i}^{j} \geq V_{ij}&\Leftrightarrow &-\mathbf{Id}_{7}(\overline{x}^{3})^{T}\leq -\mathbf{V}\vspace{0.25cm}\\
    \diamond \quad \sum_{i=1}^{4}\sum_{j=1}^{j_{i}} x_{3i}^{j} \leq \sum_{i=1}^{4} x_{2i}&\Leftrightarrow &\mathbf{1}_{11}(-\overline{x}^{2},\overline{x}^{3})^{T}\leq 0
  \end{array}
\right.
$\\
\\
 \vspace{0.5cm}
  \fbox{
  
      \textit{19 Constraints}
  }
\end{tabular}
\\
 \bottomrule
    \end{tabular}
	\vspace{0.2cm}
	\caption{Matrix Formulation of objective functions and constraints of the three levels programming Problem\label{tabFOR}} 
		\end{table}
\begin{notation}[Table \ref{tabFOR}]$\;$\\
 \begin{itemize}
    \item[$\diamond$] $\mathbf{Id}_{4}$ and $\mathbf{Id}_{7}$ are the identity matrices of order 4 and 7, respectively.
    \item[$\diamond$] $\mathbf{0}_{h_1 \times h_2}$ and $\mathbf{1}_{h_1 \times h_2}$: Matrices of order $h_1 \times h_2$, where $h_1, h_2 \in \mathbb{N}$, with all components equal to 0 and 1, respectively.
    \item[$\diamond$] Decision vector:
		$$
		x := (\overline{x}^{1}, \overline{x}^{2}, \overline{x}^{3}) = (\underbrace{x_{11}}_{\overline{x}^{1}}, \underbrace{x_{21}, \ldots , x_{24}}_{\overline{x}^{2}}, \underbrace{x_{31}^1, \ldots ,x_{31}^{8}, x_{32}^1, \ldots , x_{32}^5, x_{33}^1, \ldots , x_{33}^5, x_{34}^1, \ldots , x_{34}^4}_{\overline{x}^{3}}).
		$$
		\item[$\diamond$] Vectors: $\mathbf{P}=(P_{i})$, $\mathbf{I}=(I_{i})$, $\mathbf{C}=(C_{i})$ (all of order 4), $\mathbf{B}=(B_{ij})$, $\mathbf{V}=(V_{ij})$ (both of order 7).

  \end{itemize}
\end{notation}

	\subsection{Numerical Application \label{MDLFA}}

Let's consider the three-level linear programming problem of the Vaccination Planning in Long-Term Care Facilities as follows:

\begin{equation}\label{exemple13231}
\begin{array}{ll}
 \textbf{Level}\hspace{0.12cm} 1 \textbf{:}& {\color{blue}\underset{\color{orange}\overline{x}^{1}}{\Max}}\hspace{0.1cm} f_{1}(x)=f_{1}\left({\color{orange}\overline{x}^{1}},\overline{x}^{2},\overline{x}^{3}\right)={\color{orange}x_{11}}\\
& \text{such that}\hspace{0.1cm}x_{21}, \ldots , x_{24}, \text{and}\hspace{0.1cm} x_{31}^1, \ldots ,x_{31}^{8}, x_{32}^1, \ldots , x_{32}^5, x_{33}^1, \ldots , x_{33}^5, x_{34}^1, \ldots , x_{34}^4\hspace{0.1cm} \text{solves:}\\
\textbf{Level} \hspace{0.12cm}2 \textbf{:} & {\color{blue}\underset{\color{orange}\overline{x}^{2}}{\Max}}\hspace{0.1cm} f_{2}(x)=f_{2}\left(\overline{x}^{1},{\color{orange}\overline{x}^{2}},\overline{x}^{3}\right)=\lambda_{1}{\color{orange}x_{21}}+\lambda_{2}{\color{orange}x_{22}}+\lambda_{3}{\color{orange}x_{23}}+\lambda_{4}{\color{orange}x_{24}}\\
& \text{such that}\hspace{0.1cm} x_{31}^1, \ldots ,x_{31}^{8}, x_{32}^1, \ldots , x_{32}^5, x_{33}^1, \ldots , x_{33}^5, x_{34}^1, \ldots , x_{34}^4 \hspace{0.1cm} \text{solves:}\\
 \textbf{Level}\hspace{0.12cm} 3 \textbf{:}
&  {\color{blue}\underset{\color{orange}\overline{x}^{3}}{\Max}}\hspace{0.1cm} f_{3}(x)=f_{3}\left(\overline{x}^{1},\overline{x}^{2},{\color{orange}\overline{x}^{3}}\right)=\sum_{i=1}^{4}\sum_{j=1}^{j_{i}}\omega_{i}^{j}{\color{orange}x_{3i}^{j}},
\end{array}
\end{equation}
subject to:
$$
x\in S=\left\{x\in\mathbb{R}^{27}\hspace{0.15cm}:\hspace{0.15cm}Ax= b, x\geq 0\right\}\neq\emptyset,
$$
where $A$ is a matrix of order $29 \times 27$ and $b$ is a vector of order $29$, defined as follows:
$$
\begin{array}{cc}
A = \left(\begin{array}{crr}
        1        &  \mathbf{0}_{1\times4}        & \mathbf{0}_{1\times7}        \\
	\mathbf{0}_{4\times1} &  \mathbf{Id}_{4}       & \mathbf{0}_{4\times7} \\
	\mathbf{0}_{4\times1} & -\mathbf{Id}_{4}       & \mathbf{0}_{4\times7} \\
	\mathbf{0}_{4\times1} & -\mathbf{Id}_{4}       & \mathbf{0}_{4\times7} \\
	      -1        &  \mathbf{1}_{1\times4}        & \mathbf{0}_{1\times7}        \\
	\mathbf{0}_{7\times1} &  \mathbf{0}_{7\times4}	&  \mathbf{Id}_{7}		    \\
	\mathbf{0}_{7\times1} &  \mathbf{0}_{7\times4}	& -\mathbf{Id}_{7}		    \\
	     0         &  -\mathbf{1}_{1\times4}        & \mathbf{1}_{1\times7} 
\end{array}\right),&\hspace{0.3cm}
b=\left(\begin{array}{r}100\\\mathbf{P}\\ -\mathbf{I}\\ -\mathbf{C}\\ 0 \\\mathbf{B}\\ -\mathbf{V} \\0
\end{array}\right).
\end{array}
$$



Assume the government has selected the parameters based on the importance of care hospital capacity $\mathbf{B}$:
\begin{table}[htb]
    \centering
    \renewcommand{\arraystretch}{1} 
    \setlength{\tabcolsep}{4pt} 
    \begin{tabular}{c|c|cc|cc|cc|cc|cc|cc|cc|cc}
        \toprule
        
        \multirow{2}{*}{$\lambda_{i}$ }& \diagbox[width=3em]{}{$j$} & \multicolumn{2}{c|}{1} & \multicolumn{2}{c|}{2} & \multicolumn{2}{c|}{3} & \multicolumn{2}{c|}{4} & \multicolumn{2}{c|}{5} & \multicolumn{2}{c|}{6} & \multicolumn{2}{c|}{7} & \multicolumn{2}{c}{8} \\
\cmidrule{3-18}
        &\diagbox[width=3em]{$i$}{}& $\xi_{i}^{j}$&$\omega_{i}^{j}$&$\xi_{i}^{j}$&$\omega_{i}^{j}$&$\xi_{i}^{j}$&$\omega_{i}^{j}$&$\xi_{i}^{j}$&$\omega_{i}^{j}$&$\xi_{i}^{j}$&$\omega_{i}^{j}$&$\xi_{i}^{j}$&$\omega_{i}^{j}$&$\xi_{i}^{j}$&$\omega_{i}^{j}$&$\xi_{i}^{j}$&$\omega_{i}^{j}$\\
        \midrule
        0.4 & 1 & \num{0.12}&\num{0.05} & \num{0.14}&\num{0.06} & \num{0.10}&\num{0.04} & \num{0.25}&\num{0.08} & \num{0.13}&\num{0.09} & \num{0.12}&\num{0.05} & \num{0.15}&\num{0.06} & \num{0.22}&\num{0.07} \\
        \addlinespace
        0.3 & 2 & \num{0.25}&\num{0.08} & \num{0.16}&\num{0.10} & \num{0.21}&\num{0.13} & \num{0.18}&\num{0.06} & \num{0.20}&\num{0.12}  \\
        \addlinespace
        0.2 & 3 & \num{0.31}&\num{0.15} & \num{0.23}&\num{0.12} & \num{0.18}&\num{0.09} & \num{0.15}&\num{0.06} & \num{0.25}&\num{0.08} \\
        \addlinespace
        0.1 & 4 & \num{0.30}&\num{0.03} & \num{0.27}&\num{0.02} & \num{0.23}&\num{0.07} & \num{0.20}&\num{0.02} \\
        \bottomrule
    \end{tabular}
		\vspace{0.2cm}
    \caption{Coefficients of objective functions $f_{2}$ and $f_{3}$ (scaled by $10^{4}$)}
\end{table}

Using the monocriteria simplex method, we maximize each of the three objective functions, $f_1$, $f_2$, and $f_3$. The results of this optimization process are presented in Table \ref{tab4}.

\begin{table}[htb]
  
   \renewcommand{\arraystretch}{1.25} 
    \setlength{\tabcolsep}{5.6pt} 
    	 \begin{tabular}{|c|c|c|cccc|c|c|c|c|c|c|c|c|c}
        \toprule
        {\bf DMs} & $f_{p}^{\text{max}}$ & $\overset{o}{x}_{11}$&$\overset{o}{x}_{21}$&$\overset{o}{x}_{22}$&$\overset{o}{x}_{23}$&$\overset{o}{x}_{24}$&$\overset{o}{x}_{31}^{1}$&$\overset{o}{x}_{31}^{2}$&$\overset{o}{x}_{31}^{3}$&$\overset{o}{x}_{31}^{4}$&$\overset{o}{x}_{31}^{5}$&$\overset{o}{x}_{31}^{6}$&$\overset{o}{x}_{31}^{7}$&$\overset{o}{x}_{31}^{8}$\\
        \midrule
        1& 10000 &10000 &5180&524&302&183 &3&3&2&3&3&3&4&2 \\
         \midrule
                                 2  & 2505  &7289.6 &5648&545&316&191 &2.6&3.1&2.1&2.9&3.3&2.8&3.6&2.4 \\
                                  \midrule
            3& 2762.6&699.1&5147.5&524.7&302.7&183.6&3&3.5&2.5&3.2&3.4&3&3.8&2.8 \\           
        \bottomrule
    \end{tabular}
		
		  \vspace{\baselineskip}
			 \begin{tabular}{|c|c|c|c|c|c|c|c|c|c|c|c|c|c|c|}
        \toprule
        {\bf DMs} &$\overset{o}{x}_{32}^{1}$&$\overset{o}{x}_{32}^{2}$&$\overset{o}{x}_{32}^{3}$&$\overset{o}{x}_{32}^{4}$&$\overset{o}{x}_{32}^{5}$&$\overset{o}{x}_{33}^{1}$&$\overset{o}{x}_{33}^{2}$&$\overset{o}{x}_{33}^{3}$&$\overset{o}{x}_{33}^{4}$&$\overset{o}{x}_{33}^{5}$&$\overset{o}{x}_{34}^{1}$&$\overset{o}{x}_{34}^{2}$&$\overset{o}{x}_{34}^{3}$&$\overset{o}{x}_{34}^{4}$\\
        \midrule
        1&4&2&3&3&3&2&1&1&1&1&2&2&1&1 \\
        \midrule
      2&3.6&2.1&3.1&2.6&2.9&1.6&1.3&1.1  &0.9&0.7&1.8&1.6&1.4&1.2 \\
      \midrule
            3&4&2.5&3.3&2.8&3.1&2&1.5&1.2&1&0.8&1.9&1.7&1.5&1.3\\           
        \bottomrule
    \end{tabular}
		\vspace{0.2cm}
    \caption{Optimal solutions $\overset{o}{x}_{1}$, $\overset{o}{x}_{2}$, $\overset{o}{x}_{3}$ and optimal values $f_{1}^{\text{max}}$, $f_{2}^{\text{max}}$, $f_{3}^{\text{max}}$ (scaled by $10^{4}$) \label{tab4}}

\end{table}

       \begin{sidewaystable}[!htbp]
    \centering
		\renewcommand{\arraystretch}{1.25} 
    \setlength{\tabcolsep}{1.5pt} 
    
   \begin{tabular}{*{2}{c|}c|cccc|c|cccc|cccccccccccccccccccccccc}
        \toprule
\multicolumn{2}{c|}{}&\multicolumn{5}{c|}{\bf Simplex Method}&\multicolumn{5}{c|}{\bf Adaptive Method}&\multicolumn{22}{c}{\bf Adaptive method and Simplex Method}\\
        \midrule
        $\alpha_{1}$&$\alpha_{2}$&$\overset{c}{x}_{11}$&$\overset{c}{x}_{21}$&$\overset{c}{x}_{22}$&$\overset{c}{x}_{23}$&$\overset{c}{x}_{24}$&$\overset{c}{x}_{11}$&$\overset{c}{x}_{21}$&$\overset{c}{x}_{22}$&$\overset{c}{x}_{23}$&$\overset{c}{x}_{24}$&$\overset{c}{x}_{31}^{1}$&$\overset{c}{x}_{31}^{2}$&$\overset{c}{x}_{31}^{3}$&$\overset{c}{x}_{31}^{4}$&$\overset{c}{x}_{31}^{5}$&$\overset{c}{x}_{31}^{6}$&$\overset{c}{x}_{31}^{7}$&$\overset{c}{x}_{31}^{8}$&$\overset{c}{x}_{32}^{1}$&$\overset{c}{x}_{32}^{2}$&$\overset{c}{x}_{32}^{3}$&$\overset{c}{x}_{32}^{4}$&$\overset{c}{x}_{32}^{5}$&$\overset{c}{x}_{33}^{1}$&$\overset{c}{x}_{33}^{2}$&$\overset{c}{x}_{33}^{3}$&$\overset{c}{x}_{33}^{4}$&$\overset{c}{x}_{33}^{5}$&$\overset{c}{x}_{34}^{1}$&$\overset{c}{x}_{34}^{2}$&$\overset{c}{x}_{34}^{3}$&$\overset{c}{x}_{34}^{4}$  \\
        \midrule
        \multirow{5}{*}{0}&  0     & 83046 & 54423 & 5323& 3070 & 1854&69910&51475&5245&3024&1833 & 30 & 35 & 25 & 32 & 34 & 30 & 38 & 28 & 40 & 25 & 33 & 28 & 31 & 20 & 15 & 12 & 10 & 08 & 19 & 17 & 15 & 13 \\
                          &   $0.25$   & 83046 & 54423 & 5323 & 3070 & 1854 &69910&51475&5245&3024&1833& 30 & 35 & 25 & 32 & 34 & 30 & 38 & 28 & 40 & 25 & 33 & 28 & 31 & 20 & 15 & 12 & 10 & 08 & 19 & 17 & 15 & 13 \\
                          &   $0.5$   & 83046 & 54423 & 5323 & 3070 & 1854 &69910&51475&5245&3024&1833& 30 & 35 & 25 & 32 & 34 & 30 & 38 & 28 & 40 & 25 & 33 & 28 & 31 & 20 & 15 & 12 & 10 & 08 & 19 & 17 & 15 & 13 \\
                          &   $0.75$    & 83046 & 54423 & 5323 & 3070 & 1854 &69910&51475&5245&3024&1833& 30 & 35 & 25 & 32 & 34 & 30 & 38 & 28 & 40 & 25 & 33 & 28 & 31 & 20 & 15 & 12 & 10 & 08 & 19 & 17 & 15 & 13 \\\
                          &    1     & 83046 & 54423 & 5323 & 3070 & 1854 &69910&51475&5245&3024&1833& 30 & 35 & 25 & 32 & 34 & 30 & 38 & 28 & 40 & 25 & 33 & 28 & 31 & 20 & 15 & 12 & 10 & 08 & 19 & 17 & 15 & 13 \\
        \midrule
        \multirow{5}{*}{$0.25$}&0 &
				 81552 & 54376 & 5323 & 3070 & 1854 &69910&51475&5245&3024&1833& 30 & 35 & 25 & 32 & 34 & 30 & 38 & 28 & 40 & 25 & 33 & 28 & 31 & 20 & 15 & 12 & 10 & 08 & 19 & 17 & 15 & 13 \\
        
        &0.25 &81487 & 54064 & 5304 & 3057 & 1853 &69910&51475&5245&3024&1833& 30 & 35 & 25 & 32 & 34 & 30 & 38 & 28 & 40 & 25 & 33 & 28 & 31 & 20 & 15 & 12 & 10 & 08 & 19 & 17 & 15 & 13 \\
        
        &0.5 &81382 & 53502 & 5279 & 3042 & 1851 &69910&51475&5245&3024&1833& 30 & 35 & 25 & 32 & 34 & 30 & 38 & 28 & 40 & 25 & 33 & 28 & 31 & 20 & 15 & 12 & 10 & 08 & 19 & 17 & 15 & 13 \\
        
        &0.75 &81199 & 52453 & 5259 & 3033 & 1846 &69910&51475&5245&3024&1833& 30 & 35 & 25 & 32 & 34 & 30 & 38 & 28 & 40 & 25 & 33 & 28 & 31 & 20 & 15 & 12 & 10 & 08 & 19 & 17 & 15 & 13 \\
        
        &1 &81029 & 51475 & 5245 & 3024 & 1855 &69910&51475&5245&3024&1833& 30 & 35 & 25 & 32 & 34 & 30 & 38 & 28 & 40 & 25 & 33 & 28 & 31 & 20 & 15 & 12 & 10 & 08 & 19 & 17 & 15 & 13 \\
      \midrule
				
				\multirow{5}{*}{$0.5$}&0 &
				
				79557 & 54295 & 5324 & 3071 & 1854 &69910&51475&5245&3024&1833& 30 & 35 & 25 & 32 & 34 & 30 & 38 & 28 & 40 & 25 & 33 & 28 & 31 & 20 & 15 & 12 & 10 & 08 & 19 & 17 & 15 & 13 \\
        &0.25&79505 & 54004 & 5306 & 3071 & 1854 &69910&51475&5245&3024&1833& 30 & 35 & 25 & 32 & 34 & 30 & 38 & 28 & 40 & 25 & 33 & 28 & 31 & 20 & 15 & 12 & 10 & 08 & 19 & 17 & 15 & 13 \\
        &0.5&79437 & 53496 & 5285 & 3072 & 1855 &69910&51475&5245&3024&1833& 30 & 35 & 25 & 32 & 34 & 30 & 38 & 28 & 40 & 25 & 33 & 28 & 31 & 20 & 15 & 12 & 10 & 08 & 19 & 17 & 15 & 13 \\
        &0.75&79268 & 52491 & 5260 & 3060 & 1847 &69910&51475&5245&3024&1833& 30 & 35 & 25 & 32 & 34 & 30 & 38 & 28 & 40 & 25 & 33 & 28 & 31 & 20 & 15 & 12 & 10 & 08 & 19 & 17 & 15 & 13 \\
        &1&78905 & 51475 & 5245 & 3056 & 1854 &69910&51475&5245&3024&1833& 30 & 35 & 25 & 32 & 34 & 30 & 38 & 28 & 40 & 25 & 33 & 28 & 31 & 20 & 15 & 12 & 10 & 08 & 19 & 17 & 15 & 13 \\
       \midrule
\multirow{5}{*}{$0.75$}&0 &
            76778 & 54168 & 5322 & 3071 & 1854 &69910&51475&5245&3024&1833& 30 & 35 & 25 & 32 & 34 & 30 & 38 & 28 & 40 & 25 & 33 & 28 & 31 & 20 & 15 & 12 & 10 & 08 & 19 & 17 & 15 & 13 \\
       &0.25&
        76725 & 53897 & 5304 & 3058 & 1847 &69910&51475&5245&3024&1833& 30 & 35 & 25 & 32 & 34 & 30 & 38 & 28 & 40 & 25 & 33 & 28 & 31 & 20 & 15 & 12 & 10 & 08 & 19 & 17 & 15 & 13 \\
      &0.5&
        76657 & 53425 & 5283 & 3044 & 1839 &69910&51475&5245&3024&1833& 30 & 35 & 25 & 32 & 34 & 30 & 38 & 28 & 40 & 25 & 33 & 28 & 31 & 20 & 15 & 12 & 10 & 08 & 19 & 17 & 15 & 13 \\
     &0.75&
        76589 & 52501 & 5263 & 3032 & 1838 &69910&51475&5245&3024&1833& 30 & 35 & 25 & 32 & 34 & 30 & 38 & 28 & 40 & 25 & 33 & 28 & 31 & 20 & 15 & 12 & 10 & 08 & 19 & 17 & 15 & 13 \\
       &1&
        75794 & 51475 & 5245 & 3024 & 1833 &69910&51475&5245&3024&1833& 30 & 35 & 25 & 32 & 34 & 30 & 38 & 28 & 40 & 25 & 33 & 28 & 31 & 20 & 15 & 12 & 10 & 08 & 19 & 17 & 15 & 13 \\
         \midrule
       \multirow{5}{*}{$1$}&0 &
				  69910 & 54073 & 5305 & 3068 & 1863 &69910&51475&5245&3024&1833& 30 & 35 & 25 & 32 & 34 & 30 & 38 & 28 & 40 & 25 & 33 & 28 & 31 & 20 & 15 & 12 & 10 & 08 & 19 & 17 & 15 & 13 \\
        &0.25&
        69910 & 54073 & 5305 & 3068 & 1863 &69910&51475&5245&3024&1833& 30 & 35 & 25 & 32 & 34 & 30 & 38 & 28 & 40 & 25 & 33 & 28 & 31 & 20 & 15 & 12 & 10 & 08 & 19 & 17 & 15 & 13 \\
     &0.5&
        69910 & 54073 & 5305 & 3068 & 1863 &69910&51475&5245&3024&1833& 30 & 35 & 25 & 32 & 34 & 30 & 38 & 28 & 40 & 25 & 33 & 28 & 31 & 20 & 15 & 12 & 10 & 08 & 19 & 17 & 15 & 13 \\
      &0.75&
        69910 & 54073 & 5305 & 3068 & 1863 &69910&51475&5245&3024&1833& 30 & 35 & 25 & 32 & 34 & 30 & 38 & 28 & 40 & 25 & 33 & 28 & 31 & 20 & 15 & 12 & 10 & 08 & 19 & 17 & 15 & 13 \\
       &1&
        69910 & 54073 & 5305 & 3068 & 1863 &69910&51475&5245&3024&1833& 30 & 35 & 25 & 32 & 34 & 30 & 38 & 28 & 40 & 25 & 33 & 28 & 31 & 20 & 15 & 12 & 10 & 08 & 19 & 17 & 15 & 13 \\

        \bottomrule
\end{tabular}

    \caption{Comprehensive overview of solutions for the Linear three-level problem \ref{exemple13231} obtained through simplex method and adapted method across various $\alpha_{1}$ and $\alpha_{2}$ values (scaled by $10^{3}$)\label{tab5}}
\end{sidewaystable}
Applying our algorithm \ref{alg:mlmolpp}, given that we have three $DMs$ (3 levels, $P=3$), we will undergo two iterations of the while loop in step 3. The $DM_{1}$ and $DM_{2}$ will, in turn, narrow down the intervals of the components under their control. The compromise will be reached by $DM_{3}$, and the algorithm terminates. Let $\alpha_{1}$ and $\alpha_{2}$ denote the parameters chosen by the first two $DMs$, which reduce their respective intervals, specifically in the case of $\overset{o}{x}_{p-1,j}^{p-1}$. It is important to note that both $\alpha8{1}$ and $\alpha_{2}$ are constrained to values between 0 and 1
\[
\begin{array}{cc}
\left[l_{p,j}^{p},u_{p,j}^{p}\right] =
\begin{cases}
\left[l_{p-1,j}^{p-1}+\alpha_{p-1}\left|u_{p-1,j}^{p-1}-l_{p-1,j}^{p-1}\right|,u_{p-1,j}^{p-1}\right], & \text{if } l_{p-1,j}^{p-1}=\overset{o}{x}_{p-1,j}^{p-1},\vspace{0.2cm} \\
\left[l_{p-1,j}^{p-1},u_{p-1,j}^{p-1}-\alpha_{p-1}\left|u_{p-1,j}^{p-1}-l_{p-1,j}^{p-1}\right|\right], & \text{if } u_{p-1,j}^{p-1}=\overset{o}{x}_{p-1,j}^{p-1},
\end{cases}
& \begin{array}{l}p=2,3, \\ j=1,\ldots,n_{p-1}.\end{array}
\end{array}
\]
The compromise solutions $\overset{c}{x}$ are then obtained as functions of $\alpha_{1}$ and $\alpha_{2}$, as detailed in Table \ref{tab5}.

\subsubsection{Table \ref{tab6} Explanation: Comparative Analysis of Simplex and Adaptive Methods}		
Within the specified range of [0,1], the parameters $\alpha_{1}$ and $\alpha_{2}$ are carefully chosen with a subdivision step of [0,1], set at 0.25. Subsequently, a comprehensive exploration of all possible combinations of $\alpha_{1}$ and $\alpha_{2}$ is conducted by the two $DM$s, leading to the data presented in Table \ref{tab6}. The performance metrics under scrutiny include time efficiency, memory consumption, and objective function values. The computational resolution is executed using Matlab 2007b on a PC featuring an 'Intel(R) Celeron(R) CPU N3150 @ 1.60GHz 1.60 GHz' processor with 4 GB of RAM. Time is measured in seconds, while memory is recorded in megabytes (Mb). Notably, the objective function values ($\overset{c}{f}_1$, $\overset{c}{f}_2$, $\overset{c}{f}_3$) are scaled by a factor of $10^{4}$ for improved readability. These objective function values signify the maximum outcomes attained by each method concerning the respective objective functions, computed over the compromise solutions $\overset{c}{x}$.



\begin{table}[htb]
    \centering
    	\renewcommand{\arraystretch}{1} 
    \setlength{\tabcolsep}{5.5pt} 
    
    \begin{tabular}{c|c|ccccc|ccccc}
         \toprule
				\multicolumn{2}{c|}{}&\multicolumn{5}{c|}{{\bf Simplex method}}&\multicolumn{5}{c}{{\bf Adaptive method}}\\
      
     \midrule
					        $\alpha_{1}$    &  $\alpha_{2}$          & {\bf Time} &{\bf Memory }& \multicolumn{3}{c|}{{\bf Values} ($\times 10^4$)}&  {\bf Time} &{\bf Memory}&  \multicolumn{3}{c}{{\bf Values} ($\times 10^4$)}\\
											&						 &		     (s)       &     (Mb)         &$\overset{c}{f}_{1}$&$\overset{c}{f}_{2}$&$\overset{c}{f}_{3}$&       (s)       &        (Mb)      &$\overset{c}{f}_{1}$&$\overset{c}{f}_{2}$& $\overset{c}{f}_{3}$\\
				 \midrule
        \multirow{5}{*}{0} & 0   &   0.9751   & 105572 & 8304.6  & 2416.5 &2.7626&   1.2137   & 96002  &6991&2295.1&2.7626 \\
                           &0.25 &   0.9491   & 49298  & 8304.6  & 2416.5 &2.7626&   0.9276   & 59340  &6991&2295.1&2.7626 \\  
													 & 0.5 &   0.9507   & 16576  & 8304.6  & 2416.5 &2.7626&   0.9162   & 16576  &6991&2295.1&2.7626 \\
													 &0.75 &   0.9564   & 16576  & 8304.6  & 2416.5 &2.7626&   0.9139   & 16576  &6991&2295.1&2.7626 \\
												   &  1  &   0.9639   & 16576  & 8304.6  & 2416.5 &2.7626&   0.9215   & 16576  &6991&2295.1&2.7626 \\
        \midrule
     \multirow{5}{*}{0.25} &  0   &  0.9534   & 16616  & 8155.2  & 2416.5 &2.7626&   0.9206   & 16616  &6991&2295.1&2.7626 \\
                           &0.25  &  0.9600   & 16576  & 8148.7  & 2401.3 &2.7626&   0.9177   & 16576  &6991&2295.1&2.7626 \\ 
													 & 0.5  &  0.9596   & 16576  & 8138.2  & 2377.8 &2.7626&   0.9186   & 16576  &6991&2295.1&2.7626 \\
													 &0.75	&  0.9532   & 16576  & 8119.9  & 2335.0 &2.7626&   0.9146   & 16576  &6991&2295.1&2.7626 \\
													 &  1   &  0.9658   & 16576  & 8102.9  & 2295.4 &2.7626&   0.9160   & 16576  &6991&2295.1&2.7626 \\
        \midrule
		  \multirow{5}{*}{0.5} &  0   &  0.9500   & 16616  & 7955.7  & 2411.5 &2.7626&   0.9201   & 16616  &6991&2295.1&2.7626 \\
                           &0.25  &  0.9553   & 16576  & 7950.5  & 2399.3 &2.7626&   0.9175   & 16576  &6991&2295.1&2.7626 \\  
													 & 0.5  &  0.9582   & 16576  & 7943.7  & 2378.3 &2.7626&   0.9205   & 16576  &6991&2295.1&2.7626 \\
													 &0.75	&  0.9597   & 16576  & 7926.8  & 2337.1 &2.7626&   0.9194   & 16576  &6991&2295.1&2.7626 \\
												   &  1   &  0.9641   & 16576  & 7890.5  & 2296.0 &2.7626&   0.9185   & 16576  &6991&2295.1&2.7626 \\
        \midrule
		 \multirow{5}{*}{0.75} &  0   &  0.9581   & 16616  & 7677.8  & 2406.3 &2.7626&  0.9212   & 16616  &6991&2295.1&2.7626 \\
                           &0.25  &  0.9604   & 16576  & 7672.5  & 2394.6 &2.7626&  0.9212   & 16576  &6991&2295.1&2.7626 \\ 
													 & 0.5  &  0.9532   & 16576  & 7665.7  & 2374.8 &2.7626&  0.9256   & 16576  &6991&2295.1&2.7626 \\
													 &0.75  &  0.9592   & 16576  & 7658.9  & 2337.0 &2.7626&  0.9178   & 16576  &6991&2295.1&2.7626 \\
													 &  1   &  0.9825   & 16576  & 7579.4  & 2295.1 &2.7626&  0.9296   & 16576  &6991&2295.1&2.7626 \\
        \midrule
				\multirow{5}{*}{1} &  0   &  0.9534   & 16616  & 6991    & 2402.1 &2.7626&  0.9222   & 16616  &6991&22951&2.7626 \\
                           &0.25  &  0.9603   & 16576  & 6991    & 2402.1 &2.7626&  0.9158   & 16576  &6991&22951&2.7626 \\
													 & 0.5  &  0.9556   & 16576  & 6991    & 2402.1 &2.7626&  0.9191   & 16576  &6991&22951&2.7626 \\
													 &0.75	&  0.9688   & 16576  & 6991    & 2402.1 &2.7626&  0.9317   & 16576  &6991&22951&2.7626 \\
													 &  1   &  0.9601   & 16576  & 6991    & 2402.1 &2.7626&  0.9169   & 16576  &6991&22951&2.7626 \\

        \bottomrule
    \end{tabular}
		\vspace{0.2cm}
   \caption{Results of the Simplex and Adaptive methods for varying $\alpha_{1}$ and $\alpha_{2}$. The values in the "Values" columns are scaled by $10^4$.\label{tab6}}
\end{table}

\begin{figure}[htbp]
{ 
    \centering

\begin{minipage}[b]{0.37\textwidth}
        \begin{subfigure}[b]{\textwidth}
            \includegraphics[width=6cm,height=3.15cm]{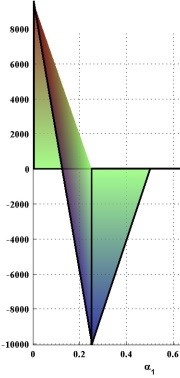}
            \caption{Memory Allocation at $\alpha_2=0$}
        \end{subfigure}
        \vfill 
        \begin{subfigure}[b]{\textwidth}
            \includegraphics[width=6cm,height=3.15cm]{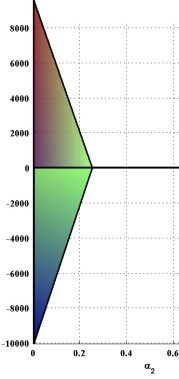}
            \caption{Memory Allocation at $\alpha_1=0$}
        \end{subfigure}
    \end{minipage} \hfill
   \begin{minipage}[b]{0.6\textwidth}
        \begin{subfigure}[b]{\textwidth}
            \includegraphics[width=9cm,height=7.5cm]{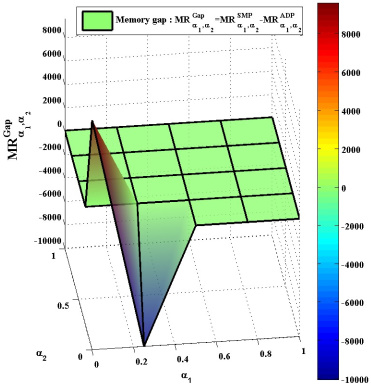}
            \caption{Spatial Insights into Memory Allocation: Contrasting Simplex and Adapted Method}
        \end{subfigure}
    \end{minipage}
    \caption{Memory Allocation Variability: A Comparative Study between Simplex and Adapted Methods, Influenced by $\alpha_1$ and $\alpha_2$}
    \label{fig:memo}
}

 \paragraph{Memory Utilization Discrepancy:} The variable $MR_{\alpha_{1},\alpha_{2}}^{\text{Gap}}$ characterizes the difference in memory usage between the Adaptive (ADP) method  and the Simplex (SMP) method . In this context, $MR_{\alpha_{1},\alpha_{2}}^{\text{ADP}}$ represents the memory consumption of the ADP method under parameters $\alpha_{1}$ and $\alpha_{2}$, while $MR_{\alpha_{1},\alpha_{2}}^{\text{SMP}}$ represents the memory consumption of the SMP method under the same parameters.
\end{figure}


\begin{figure}[htbp]
{ 
    \centering

\begin{minipage}[b]{0.37\textwidth}
        
            \includegraphics[width=5cm,height=1cm]{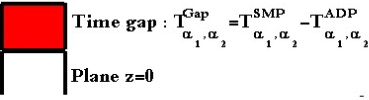}
       
        \vfill 
        \begin{subfigure}[b]{\textwidth}
            \includegraphics[width=6cm,height=5cm]{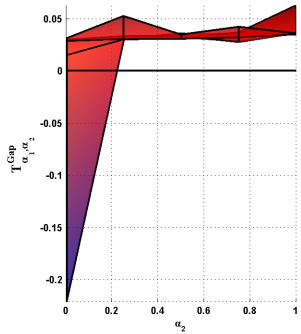}
            \caption{Planar Insights into Simplex vs. Adaptive Method Dynamics at $\alpha_{1}=0$}
        \end{subfigure}
    \end{minipage} \hfill
   \begin{minipage}[b]{0.6\textwidth}
        \begin{subfigure}[b]{\textwidth}
            \includegraphics[width=9cm,height=7cm]{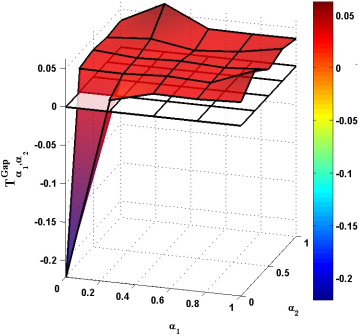}
            \caption{Temporal Efficiency Analysis: A 3D Exploration of the Simplex Method vs. an Adaptive Approach}
        \end{subfigure}
    \end{minipage}
    \caption{Comparison of Temporal Discrepancies Between the Simplex Method and an Adaptive Method Based on $\alpha_1$ and $\alpha_2$ Parameters}
    \label{fig:tmp}
    }

    \paragraph{Temporal Disparity Analysis:} The variable $T_{\alpha_{1},\alpha_{2}}^{\text{Gap}}$ delineates the temporal deviation between the ADP and the SMP methods. Here, $T_{\alpha_{1},\alpha_{2}}^{\text{ADP}}$ represents the temporal consumption of the ADP method under parameters $\alpha_{1}$ and $\alpha_{2}$, while $T_{\alpha_{1},\alpha_{2}}^{\text{SMP}}$ corresponds to the temporal consumption of the SMP method under the same parameters.
\end{figure}

\begin{figure}[htbp]
{ 
    \centering
    
\begin{minipage}[b]{0.37\textwidth}
        
            \includegraphics[width=5cm,height=1cm]{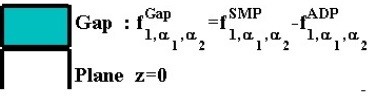}
       
        \vfill 
        \begin{subfigure}[b]{\textwidth}
            \includegraphics[width=7cm,height=4cm]{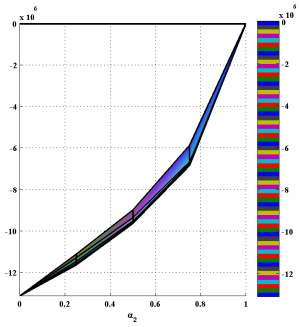}
            \caption{Snapshot at $\alpha_1=0$: Maximizing $f_1$ - Simplex vs. Adaptive Method}
        \end{subfigure}
    \end{minipage} \hfill
   \begin{minipage}[b]{0.6\textwidth}
        \begin{subfigure}[b]{\textwidth}
            \includegraphics[width=9cm,height=6cm]{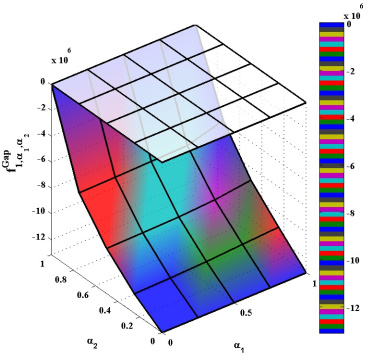}
            \caption{Solution Peaks for $f_1$: Simplex vs. Adaptive Method (Variations in $\alpha_1$ and $\alpha_2$)}
        \end{subfigure}
    \end{minipage}
    \vfill
    \begin{minipage}[b]{0.37\textwidth}
     \includegraphics[width=5cm,height=1cm]{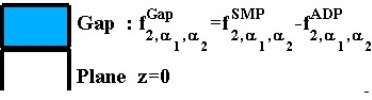}
       
        \vfill 
        
        \begin{subfigure}[b]{\textwidth}
            \includegraphics[width=7cm,height=4cm]{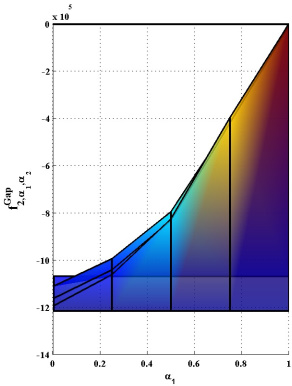}
            \caption{Visualizing $\alpha_2=0$: Peak Values in $f_2$ -  Simplex vs. Adaptive Method (Variations in $\alpha_1$ and $\alpha_2$)}
        \end{subfigure}
    \end{minipage} \hfill
   \begin{minipage}[b]{0.6\textwidth}
        \begin{subfigure}[b]{\textwidth}
            \includegraphics[width=9cm,height=6cm]{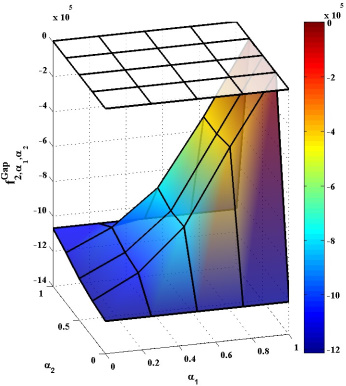}
            \caption{Maximums in $f_2$: Comparative Study between Simplex and Adaptive Method (Different $\alpha_1$ and $\alpha_2$)}
        \end{subfigure}
    \end{minipage}

     \vfill
    \begin{minipage}[b]{0.37\textwidth}
     \includegraphics[width=5cm,height=1cm]{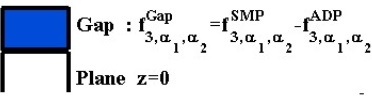}
       
        \vfill 
        
        \begin{subfigure}[b]{\textwidth}
            \includegraphics[width=7cm,height=4cm]{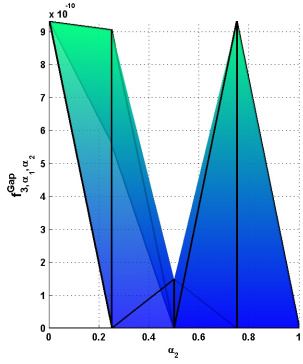}
            \caption{Snapshot at $\alpha_1=0$: Maximizing $f_3$ - Simplex vs. Adaptive Method}
        \end{subfigure}
    \end{minipage} \hfill
   \begin{minipage}[b]{0.6\textwidth}
        \begin{subfigure}[b]{\textwidth}
            \includegraphics[width=9cm,height=6cm]{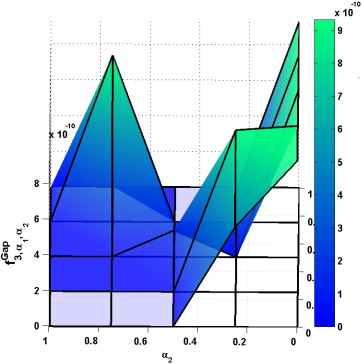}
            \caption{Optimal Peaks in $f_3$: Analyzing the Difference between Simplex and Adaptive Method (Varied $\alpha_1$ and $\alpha_2$)}
        \end{subfigure}
    \end{minipage}
    \caption{Divergence in Objective Function Values: Comparative Analysis of Simplex and Adaptive Methods with Respect to $\alpha_1$ and $\alpha_2$}
    \label{fig:div}
    }
    \paragraph{Objective Function Disparity:} The variable $f_{p,\alpha_{1},\alpha_{2}}^{\text{Gap}}$ signifies the difference in the maximum values achieved using the SMP method and the ADP method, where $p = 1, 2, 3$. Here, $f_{p,\alpha_{1},\alpha_{2}}^{\text{ADP}}$ and $f_{p,\alpha_{1},\alpha_{2}}^{\text{SMP}}$ respectively denote the maximum values attained by the ADP and SMP methods' solutions under parameters $\alpha_{1}$ and $\alpha_{2}$.
\end{figure}

\subsection{Discussion}

\subsubsection{Hospitals and Care Capacity}

As per information from the NHS website, the UK boasts around 2,500 hospitals. The data presented in Table \ref{tab:real-world-data} corresponds to 22 of these hospitals, constituting approximately 0.88\% of the total. It's crucial to note that the table's coverage is restricted to hospitals with available data regarding vulnerable populations and care capacity. Consequently, the actual count of hospital beds in the UK might surpass the provided figure, and the 22 hospitals may represent only a fraction of the total care capacity. As a result, pinpointing the precise percentage of the UK's total care capacity represented by the hospitals in the table remains challenging, but it is likely to be less than 0.88\%.

\subsubsection{Vaccine Planning in LTCFs}

The mathematical model outlined in Section \ref{MDLF}, employed in the vaccine planning example, offers a simplified portrayal of the real-world problem. Significantly, it does not encompass all factors that influence vaccine distribution, such as the prevailing epidemiological situation, vaccination policies, and resource availability. The model serves an illustrative purpose, demonstrating the application of the method to vaccine planning in LTCFs. Its goal is to showcase how the method could be applied to devise effective and equitable distribution plans in real-world scenarios, recognizing the need for further refinement to account for additional complexities in practice.

To refine the model, various enhancements can be contemplated, such as updating it to align with the present epidemiological landscape, adjusting it to accommodate evolving vaccination policies, and broadening its scope to encompass additional variables influencing vaccine distribution. While these enhancements promise a more precise depiction of the intricate facets of vaccine planning, it's essential to acknowledge that they may result in a non-linear model. Consequently, this non-linearity poses a challenge in applying our proposed method in such cases.
\subsubsection{Discussion on Optimization Results}

The results presented in Table~\ref{tab6} illuminate the comparative performance of the simplex and adaptive methods under varying parameter configurations ($\alpha_{1}$ and $\alpha_{2}$). Noticeable trends in time and memory consumption emerge as these parameters change, revealing a consistent pattern between the two methods.

To enhance the interpretation of the results in Table~\ref{tab6}, visual aids are provided through the graphs in Figures \ref{fig:memo}, \ref{fig:tmp}, and \ref{fig:div}. These graphical representations offer a clearer and more intuitive understanding of the outcomes, contributing to a comprehensive analysis and facilitating a nuanced examination of the data.

Consideration of the computational environment is crucial in interpreting the obtained results. The utilization of Matlab 2007b on a computer with 'Intel(R) Celeron(R) CPU N3150 @ 1.60GHz 1.60 GHz' processor and 4 GB of RAM, indicative of lower-end performance and an outdated Matlab version, inevitably influences the execution time and potentially the precision of our method. Recognizing these limitations is paramount. Upgrading to a more powerful system with modern hardware specifications and an updated Matlab version is expected to significantly improve result quality, especially in terms of reduced execution time. This emphasizes the importance of considering computational resources in interpreting and generalizing our study's findings.

For the scenario where $\alpha_{1}=0$ and $\alpha_{2}=0$, representing $DM$s making no reduction to their intervals (maximizing flexibility for followers), the adaptive method exhibits a longer execution time but requires less memory compared to the simplex method. However, for other values of $\alpha_{1}$ (ranging from 0 to 1) and $\alpha_{2}$ (ranging from 0.25 to 1), the adaptive method consistently demonstrates reduced time consumption compared to the simplex method. Regarding memory usage, the adaptive method consumes more memory than the simplex method only for $\alpha_{1}=0$ and $\alpha_{2}=0.25$. Subsequently, both methods consume the same amount of memory for the remaining values. This observation is reflected in the peaks in the memory graph \ref{fig:memo} and the negative spike in the time graph  \ref{fig:tmp}, followed by stabilization at positive values, indicating an advantage in favor of the adaptive method. These results underscore the significance of parameter configuration in determining the relative performance of the two methods.

An intriguing perspective is to explain why this influence particularly manifests itself in the neighborhood of the values $\alpha_{1}=\alpha_{2}=0$.

These insights are pivotal for practitioners and researchers in selecting an optimization approach based on specific requirements and resource constraints. The results contribute to a nuanced understanding of the trade-offs between the simplex and adaptive methods, guiding the selection of the most suitable method based on the optimization problem's characteristics and the available computational resources.

\subsection {Data Sources}
\begin{itemize}
\item Population data: \href{https://www.ons.gov.uk/}{Office for National Statistics (ONS)}, \href{https://www.nrscotland.gov.uk/}{National Records of Scotland (NRS)}, \href{https://gov.wales/}{Welsh Government},  \href{https://www.nisra.gov.uk/}{Northern Ireland Statistics and Research Agency (NISRA)}
\item Cases data: Public Health England, Public Health Scotland, Public Health Wales, Public Health Agency (Northern Ireland)
\item Administration capacity data: \href{https://www.nhs.uk/conditions/intensive-care/}{National Health Service (NHS)} England, Scottish Government, Welsh Government, Department of Health (Northern Ireland)
\item Equitable vaccine distribution target: \href{https://who.foundation/}{World Health Organization (WHO)}, \href{https://www.unicef.org/}{United Nations Children's Fund (UNICEF)}
\item $V_{ij}$ values were estimated by ONS, NRS, Welsh Government, and NISRA.
\item $B_{ij}$ values were provided by the respective hospitals.
\end{itemize}
\section{Conclusion}\label{section7}

This study pursued a dual objective: firstly, the application of the adaptive method to solve a ML(MO)OLPP; secondly, the refinement of Sinha and Sinha's algorithm \cite{sinha} by replacing the simplex method with the adaptive method in step 10. An additional innovation was the restructuring of steps (2)-(9) using the interval reduction map $\xi$, consolidating these eight steps into a single step (Step 3 in our algorithm).

By harnessing the constructive adaptive method, we introduced a detailed algorithm for addressing ML(MO)OLPP, with the establishment of the interval reduction map playing a pivotal role. The convergence of our algorithms is evident from the results of the adaptive method, and a numerical example was provided to illustrate the efficacy of our algorithm.

The discussions in the preceding sections have illuminated the practical implications of our work. While our algorithm exhibits enhanced efficiency in solving ML(MO)OLPP, it is essential to recognize both its scope and limitations. The comparisons in Table~\ref{tab6} offer informations on the computational performance of the simplex and adaptive methods under different parameter configurations ($\alpha_{1}$ and $\alpha_{2}$). These insights serve as a guide for selecting the most suitable optimization approach based on specific requirements and resource constraints.

In conclusion, our contributions extend beyond algorithmic enhancements. We have not only presented a comprehensive understanding of the trade-offs between optimization methods but also discussed their implications for real-world problem-solving, exemplified through a numerical example involving the distribution of COVID-19 vaccines in long-term care facilities. This application, specifically in healthcare planning, not only advances the field of optimization but also encourages future research endeavors that consider the broader context of application.

\bibliographystyle{unsrt}  
\bibliography{references}

\begin{thebibliography}{10}

\bibitem{kaci}
M.~{Kaci and S. Radjef}.
\newblock The set of all the possible compromises of a multi-level
  multi-objective linear programming problem.
\newblock {\em Croatian Operational Research Review}, 13(1):13--30, 2022.

\bibitem{kaci2023}
M.~{Kaci and S. Radjef}.
\newblock An adaptive method to solve multilevel multiobjective linear
  programming problems.
\newblock {\em Operations Research and Decisions}, 33(3):29--44, 2023.

\bibitem{abo}
M.A. Abo-sinna.
\newblock Pareto optimality for bi-level programming problem with fuzzy
  parameters.
\newblock {\em Opsearch}, 38:372--393, 2001.

\bibitem{sinna2024}
Mahmoud~A. Abo-Sinna and Ibrahim~A. Baky.
\newblock Interactive balance space approach for solving multi-level
  multi-objective programming problems.
\newblock {\em Information Sciences}, 177(16):3397--3410, 2007.
\newblock Mathematical Foundation for Intelligent Technologies (InTech03).

\bibitem{baky}
I.A. {Baky}.
\newblock Solving multi-level multi-objective linear programming problems
  through fuzzy goal programming approach.
\newblock {\em Applied Mathematical Modelling}, 34(9):2377--2387, 2010.

\bibitem{l2}
D.~Mollalign, A.~Mushi, and B.~Guta.
\newblock Solving multiobjective multilevel programming problems using
  two-phase intuitionistic fuzzy goal programming method.
\newblock {\em Journal of Computational Science}, 63:1877--7503, 2022.

\bibitem{cite1}
C.O. {Pieume and P. Marcotte and L.P. Fotso and P. Siarry}.
\newblock Solving bilevel linear multiobjective programming problems.
\newblock {\em American Journal of Operations Research}, 01:214--219, 2011.

\bibitem{ref36}
J.~{Bracken and J.T. McGill}.
\newblock Mathematical programs with optimization problems in the constraints.
\newblock {\em {Oper. Res.}}, 21:37--44, 1973.

\bibitem{ref37}
J.~{Bracken and J.T. McGill}.
\newblock Technical note—a method for solving mathematical programs with
  nonlinear programs in the constraints.
\newblock {\em Operations Research}, 22(5):917--1133, 1974.

\bibitem{BEN22}
Omar Ben-Ayed.
\newblock Bilevel linear programming.
\newblock {\em Computers and Operations Research}, 20(5):485--501, 1993.

\bibitem{BEN39}
B.~{Colson and P. Marcotte and G. Savard }.
\newblock Bilevel programming: A survey.
\newblock {\em 4OR}, 3:87–107, 2005.

\bibitem{BEN89}
V.V. {Kalashnikov, S. Dempe, G.A. Pérez-Valdés, N.I. Kalashnykova and J.-F.
  Camacho-Vallejo}.
\newblock Bilevel programming and applications.
\newblock {\em Mathematical Problems in Engineering}, 2015:1--17, 2015.

\bibitem{BEN126}
M.~{Sakawa and I. Nishizaki}.
\newblock Interactive fuzzy programming for multi-level programming problems: a
  review.
\newblock {\em International Journal of Multicriteria Decision Making},
  2(3):241–266, 2012.

\bibitem{BEN154}
L.~{ Vicente, P. Calamai}.
\newblock Bilevel and multilevel programming: a bibliography review.
\newblock {\em Journal of Global Optimization}, 5:291--306, 1994.

\bibitem{baky2009}
Ibrahim~A. Baky.
\newblock Fuzzy goal programming algorithm for solving decentralized bi-level
  multi-objective programming problems.
\newblock {\em Fuzzy Sets and Systems}, 160(18):2701--2713, 2009.
\newblock Theme: Decision and Optimisation.

\bibitem{sinha2002}
S.~{Sinha}.
\newblock Fuzzy programming approach to multi-level programming problems.
\newblock {\em Fuzzy Sets and Systems}, 136(2):189--202, 2003.

\bibitem{sinha}
S.B. {Sinha and S. Sinha}.
\newblock A linear programming approach for linear multi-level programming
  problems.
\newblock {\em The Journal of the Operational Research Society},
  55(3):312--316, 2004.

\bibitem{LM1}
K.~{Lachhwani}.
\newblock On solving multi-level multi objective linear programming problems
  through fuzzy goal programming approach.
\newblock {\em OPSEARCH}, 51(4):624–637, 2014.

\bibitem{LM2}
K.~{Lachhwani}.
\newblock Solving the general fully neutrosophic multi-level multiobjective
  linear programming problems.
\newblock {\em OPSEARCH}, 58(4):1192–1216, 2021.

\bibitem{BAH}
I.A. {Baky, M.H. Eid and M.A. El Sayed}.
\newblock Bi-level multi-objective programming problem with fuzzy demands: a
  fuzzy goal programming algorithm.
\newblock {\em OPSEARCH}, 51(2):280–296, 2014.

\bibitem{cov1}
Nezir {Aydin and Zeynep Cetinkale}.
\newblock Analyses on icu and non-icu capacity of government hospitals during
  the covid-19 outbreak via multi-objective linear programming: An evidence
  from istanbul.
\newblock {\em Computers in Biology and Medicine}, 146:105562, 2022.

\bibitem{cov2}
Sujata {Pardeshi, Sushopti Gawade, Palivela Hemant}.
\newblock Student learning time analysis during covid-19 using linear
  programming - simplex method.
\newblock {\em Social Sciences and Humanities Open}, 5(1):100266, 2022.

\bibitem{cov3}
Erfan {Babaee Tirkolaee, Hêriş Golpîra, Ahvan Javanmardan and Reza Maihami}.
\newblock A socio-economic optimization model for blood supply chain network
  design during the covid-19 pandemic: An interactive possibilistic programming
  approach for a real case study.
\newblock {\em Socio-Economic Planning Sciences}, 85:101439, 2023.

\bibitem{cov4}
A.~{Miniguano-Trujillo, F. Salazar, R. Torres et al}.
\newblock An integer programming model to assign patients based on mental
  health impact for tele-psychotherapy intervention during the covid–19
  emergency.
\newblock {\em Health Care Manag Sci}, 24(2):286–304, 2021.

\bibitem{cov5}
M.~{Tavana, K. Govindan, A.K. Nasr et al}.
\newblock A mathematical programming approach for equitable covid-19 vaccine
  distribution in developing countries.
\newblock {\em Ann Oper Res}, 2021.

\bibitem{Gabasov}
R.~{Gabasov and F.M. Kirillova and S.V. Prischepova}.
\newblock {\em Optimal Feedback Control}.
\newblock Lecture Notes in Control and Information Sciences. Springer, Berlin,
  Heidelberg, 1995.

\bibitem{GABASOV1979}
R.~{Gabasov and F.M. Kirillova}.
\newblock New linear programming methods and their application to optimal
  control.
\newblock {\em IFAC Proceedings Volumes}, 12(2):17--30, 1979.
\newblock IFAC Workshop on Control Applications of Nonlinear Programming,
  Denver, USA, 21 June.

\bibitem{Gabasov2}
R.~{Gabasov and F.M. Kirillova}.
\newblock {\em Linear Programming methods}.
\newblock House, Minsk, 1978.

\end{thebibliography}
\end{document}